\documentclass{scrartcl}
\usepackage{amsmath}
\usepackage{amsfonts}
\usepackage{amssymb}
\usepackage{dsfont}
\usepackage[T1]{fontenc}
\usepackage[latin1]{inputenc}
\usepackage{epsfig}
\usepackage{graphicx}

\usepackage[boxed, lined]{algorithm2e}
\usepackage{float}
\usepackage{comment}

\newcommand{\ba}{\mathbf{a}}
\newcommand{\bb}{\mathbf{b}}
\newcommand{\bc}{\mathbf{c}}

\newcommand{\D}{{\mathcal{D}}}

\newcommand{\N}{\mathbb{N}}

\newcommand{\R}{\mathbb{R}}

\newcommand{\Rd}{\mathbb{R}^d}

\newcommand{\beq}{\begin{eqnarray*}}

\newcommand{\eeq}{\end{eqnarray*}}

\newcommand{\beqm}{\begin{eqnarray}}

\newcommand{\eeqm}{\end{eqnarray}}

\newtheorem{theorem}{Theorem}
\newtheorem{corollary}{Corollary}
\newtheorem{lemma}{Lemma}

\newcommand{\EXP}{{\mathbf E}}
\newcommand{\PROB}{{\mathbf P}}

\renewcommand{\bf}{\normalfont \bfseries}
\renewcommand{\it}{\normalfont \itshape}

\allowdisplaybreaks
\begin{document}
\renewcommand{\thefootnote}{\fnsymbol{footnote}}
\newcommand{\F}{{\cal F}}
\newcommand{\Sp}{{\cal S}}
\newcommand{\G}{{\cal G}}
\newcommand{\HH}{{\cal H}}

\begin{center}

  {\LARGE \bf
  Over-parametrized deep neural networks do not generalize well
  }
\footnote{
Running title: {\it Over-parametrized neural networks}}
\vspace{0.5cm}

Michael Kohler$^{1}$
and Adam Krzy\.zak$^{2,}$\footnote{Corresponding author. Tel:
  +1-514-848-2424 ext. 3007, Fax:+1-514-848-2830}\\

{\it $^1$
Fachbereich Mathematik, Technische Universit\"at Darmstadt,
Schlossgartenstr. 7, 64289 Darmstadt, Germany,
email: kohler@mathematik.tu-darmstadt.de}

{\it $^2$ Department of Computer Science and Software Engineering, Concordia University, 1455 De Maisonneuve Blvd. West, Montreal, Quebec, Canada H3G 1M8, email: krzyzak@cs.concordia.ca}

\end{center}
\vspace{0.5cm}

\begin{center}
December 1, 2019
\end{center}
\vspace{0.5cm}

\noindent
    {\bf Abstract}\\
    Recently it was shown in several papers that backpropagation
    is able to find the global minimum of the empirical risk on
    the training data using over-parametrized
    deep neural networks. In this paper a similar result is shown
    for deep neural networks with the sigmoidal squasher activation function
    in a regression setting, and a lower bound is presented which
proves that these networks do not generalize well on a new data
in the sense that they do not achieve the optimal minimax
rate of convergence for estimation of smooth regression functions.

    \vspace*{0.2cm}


\vspace*{0.2cm}


\section{Introduction}
\label{se1}
Deep neural networks belong to the most successful approaches
in multivariate statistical applications, see, e.g.,
Schmidhuber (2015) and the literature cited therein.
Motivated by the practical success of these networks
there has been in recent years an increasing
interest in studying the corresponding estimators both practically and
theoretically. This is often done
in the context of nonparametric regression with random
design. Here,
$(X,Y)$ is an $\Rd \times \R$--valued random vector
satisfying $\EXP \{Y^2\}<\infty$, and given a sample
of $(X,Y)$ of size $n$, i.e., given a data set
\begin{equation}
  \label{se1eq1}
\D_n = \left\{
(X_1,Y_1), \dots, (X_n,Y_n)
\right\},
\end{equation}
where
$(X,Y)$, $(X_1,Y_1)$, \dots, $(X_n,Y_n)$ are i.i.d. random variables,
the aim is to construct an estimate
\[
m_n(\cdot)=m_n(\cdot, \D_n):\Rd \rightarrow \R
\]
of the regression function $m:\Rd \rightarrow \R$,
$m(x)=\EXP\{Y|X=x\}$ such that the $L_2$ error
\[
\int |m_n(x)-m(x)|^2 \PROB_X (dx)
\]
is ``small'' (see, e.g., Gy\"orfi et al. (2002)
for a systematic introduction to nonparametric regression and
a motivation for the $L_2$ error).
In  the sequel we want to use
(deep) neural networks in order to estimate a regression function.
The starting point in defining a neural network is the
choice of an activation function $\sigma: \mathbb{R} \to \mathbb{R}$.
Here, we use in the sequel
so--called squashing functions, which are nondecreasing
and satisfy $\lim_{x \rightarrow - \infty} \sigma(x)=0$
and
$\lim_{x \rightarrow  \infty} \sigma(x)=1$.
An example of a squashing function is
the so-called sigmoidal or logistic squasher
\begin{equation}
  \label{se2eq4}
\sigma(x)=\frac{1}{1+\exp(-x)} \quad (x \in \R).
\end{equation}
In applications, also unbounded activation functions are often
used, e.g., the famous ReLU activation function
\[
\sigma(x)=\max\{x,0\}.
\]
The network architecture $(L, \textbf{k})$ depends on a positive
integer $L$ called the \textit{number of hidden layers} and a
\textit{width vector} $\textbf{k} = (k_1, \dots, k_{L}) \in
\mathbb{N}^{L}$ that describes the number of neurons in the first,
second, $\dots$, $L$-th hidden layer. A multilayer feedforward neural
network with
architecture $(L, \textbf{k})$ and activation function $\sigma$ is a real-valued function
$f: \R^d \rightarrow \R$ defined by
\begin{equation}\label{se1eq2}
f_\bc(x) = \sum_{i=1}^{k_L} c_{1,i}^{(L)} \cdot f_i^{(L)}(x) + c_0^{(L)}
\end{equation}
for some $c_{1,0}^{(L)}, \dots, c_{1,k_L}^{(L)} \in \mathbb{R}$ and for $f_i^{(L)}$'s recursively defined by
\begin{equation}
  \label{se1eq3}
f_i^{(r)}(x) = \sigma\left(\sum_{j=1}^{k_{r-1}} c_{i,j}^{(r-1)} \cdot f_j^{(r-1)}(x) + c_{i,0}^{(r-1)} \right)
\end{equation}
for some $c_{i,0}^{(r-1)}, \dots, c_{i, k_{r-1}}^{(r-1)} \in \mathbb{R}$
$(r=2, \dots, L)$
and
\begin{equation}
  \label{se1eq4}
f_i^{(1)}(x) = \sigma \left(\sum_{j=1}^d c_{i,j}^{(0)} \cdot x^{(j)} + c_{i,0}^{(0)} \right)
\end{equation}
for some $c_{i,0}^{(0)}, \dots, c_{i,d}^{(0)} \in \mathbb{R}$.

In the sequel we want to use the data (\ref{se1eq1}) in order
to choose the weights $\bc=(c_{i,j}^{(s)})_{i,j,s}$
of the neural network such that the resulting
function $f_\bc$ defined by (\ref{se1eq2})--(\ref{se1eq4}) is a good
estimate of the regression function. This can be done
for instance by applying the principle of the least squares. Here
one defines a suitable class $\F_n$ of neural networks
and chooses that function from this class which minimizes the
error on the training data, i.e., one defines the so--called
least squares neural network estimate by
\[
m_n(\cdot)= \arg \min_{f \in \F_n}\frac{1}{n} \sum_{i=1}^n |f(X_i)-Y_i|^2.
\]
Recently it was shown in several articles, that such least squares
estimates based on deep neural networks achieve  nice rates
of convergence if suitable structural constraints
on the regression function are imposed, cf., e.g.,
Kohler and Krzy\.zak (2017), Bauer and Kohler (2019),
Kohler and Langer (2019) and
Schmidt-Hieber (2019).
Eckle and Schmidt-Hieber (2019) and
Kohler, Krzy\.zak and Langer (2019)
showed that the least squares neural network
regression estimates based on deep neural networks
can achieve the rate of convergence results similar to piecewise
polynomial partition estimates where the partition is chosen
in an optimal way.
Results concerning  estimation by neural networks of piecewise polynomial regression functions with partitions having rather general smooth boundaries have been
obtained by Imaizumi and Fukamizu (2019).

Unfortunately it is not possible to compute the least squares
neural networks regression estimate exactly,  because
such computation requires minimization of the non-convex
and nonlinear function
\[
F_n(\bc) = \frac{1}{n} \sum_{i=1}^n |f_\bc(X_i)-Y_i|^2
\]
with respect to the weight vector $\bc$. In practice,
one uses gradient descent in order to compute the minimum
of the above function approximately. Here one chooses a
random starting value $\bc^{(0)}$ for the weight
vector, and then defines
\[
\bc^{(t+1)}=\bc^{(t)} - \lambda_n \cdot (\nabla_\bc F_n)(\bc^{(t)})
\quad (t=0, \dots, t_n-1)
\]
for some suitably chosen stepsize $\lambda_n>0$ and the
number of gradient descent steps $t_n \in \N$. Then
the regression estimate is defined by
$m_n(x)=f_{\bc^{(t_n)}}(x)$.

There are quite a few papers which try to prove that
backpropagation works theoretically for deep neural
networks.
The most popular approach in this context is the so--called
landscape approach.
Choromanska et al. (2015)
used random matrix theory to
derive a heuristic argument showing
that the  risk  of  most  of  the  local  minima  of the
empirical $L_2$ risk $F_n(\bc)$ is  not  much
larger  than  the  risk  of  the  global  minimum. For neural
networks with special activation function it was possible to
validate this claim, see, e.g., Arora et al. (2018),
Kawaguchi (2016),
and Du and Lee (2018), which have
analyzed gradient descent for neural networks with linear or quadratic
activation function.
But for such neural networks there do not exist good approximation results,
consequently, one cannot derive from these results
good rates of convergence for neural network regression estimates.
Du et al. (2018)
analyzed gradient descent applied to neural networks
with one hidden layer in case of an input with a Gaussian distribution.
They used the expected gradient instead of the gradient
in their gradient descent routine,
and therefore,
their result cannot be used to derive the rate of convergence results
for a neural network regression estimate
learned by the gradient descent.
Liang et al. (2018)
applied gradient descent to a modified loss function in classification,
where it is assumed that the data can be interpolated by a neural network.
Here, as we will show in this paper (cf., Theorem \ref{th2} below),
the last assumption does not lead to good rates of convergence
in nonparametric regression,
and it is unclear whether the main idea (of simplifying the estimation by a modification of the loss function)
can also be used in a regression setting.

Recently it was shown in several papers, see, e.g.,
Allen-Zhu, Li and Song (2019), Kawaguchi and Huang (2019)
and the literature cited therein, that gradient descent leads to
a small empirical $L_2$ risk in over-parametrized neural networks.
Here the results in Allen-Zhu, Li and Song (2019) are proven for
the ReLU activation function and neural networks with a polynomial
size in the sample size. The neural networks in
Kawaguchi and Huang (2019) use squashing activation functions and
are much smaller (in fact, they require only a
linear size in the sample size). In contrast to
Allen-Zhu, Li and Song (2019)
there the learning rate is set to zero for all neurons except
for neurons in the output layer, so actually they compute a linear
least squares estimate with gradient descent, which is not used
in practice.

In this paper we show a related result for a deep neural network
regression estimate with the logistic squasher activation  function, where
the learning rate is nonzero for all neurons of the network.
By analyzing the minimax rate of convergence of this estimate
in case of a general design distribution
we are able to show that this estimate does not generalize well
to new (independent) data in a sense that it does not achieve the optimal
minimax rate of convergence in case of a smooth regression function.
Here the main trick is that we also allow discrete design
distributions and prove a general result which shows that any
estimate which achieves with high probability a very small error
on  the training data in case of such distributions 
does not achieve the optimal minimax error. This is in contrast
to a recent trend in machine learning, where one tries to argue
that such estimates can achieve good rates of convergence (see,
e.g., Bartlett et al. (2019),
Belkin et al. (2019),
Hastie et al. (2019)
and the literature cited therein). We would like to point out
that our result above is not a contradiction to Belkin, Rakhlin and Tsybakov
(2018), who show that learning method which interpolates the
training data can achieve the optimal rates for nonparametric regression
problems, because it is   assumed there that the design variable
has a density with respect to the Lebesgue-Borel measure, which
is bounded away from zero and infinity.

Throughout the paper, the following notation is used:
The sets of natural numbers, natural numbers including $0$,
and real numbers
are denoted by $\N$, $\N_0$ and $\R$, respectively.
The  Euclidean norm of $x \in \Rd$
is denoted by $\|x\|$ and $\|x\|_\infty$
denotes its supremum norm.
For $f:\R^d \rightarrow \R$
\[
\|f\|_\infty = \sup_{x \in \R^d} |f(x)|
\]
is its supremum norm. 
Let $p=q+s$ for some $q \in \N_0$ and $0< s \leq 1$.
A function $f:\R^d \rightarrow \R$ is called
$(p,C)$-smooth, if for every $\alpha=(\alpha_1, \dots, \alpha_d) \in
\N_0^d$
with $\sum_{j=1}^d \alpha_j = q$ the partial derivative
$\frac{
\partial^q f
}{
\partial x_1^{\alpha_1}
\dots
\partial x_d^{\alpha_d}
}$
exists and satisfies
\[
\left|
\frac{
\partial^q f
}{
\partial x_1^{\alpha_1}
\dots
\partial x_d^{\alpha_d}
}
(x)
-
\frac{
\partial^q f
}{
\partial x_1^{\alpha_1}
\dots
\partial x_d^{\alpha_d}
}
(z)
\right|
\leq
C
\cdot
\| x-z \|^s
\]
for all $x,z \in \R^d$.

The outline of this paper is as follows: In Section \ref{se2} the
over-parametrized neural network regression estimates are defined.
The main results are presented in Section \ref{se3} and proven
in Section \ref{se4}.

\section{An over-parametrized neural network regression estimator}
\label{se2}
In the sequel we use the logistic squasher $\sigma(x)=1/(1+e^{-x})$
as the activation function, and we use a network topology where
we compute the linear combination of $k_n$ fully connected
neural networks with $L$ layers and $k_0$ neurons per layer.
Thus we define our neural networks by
\begin{equation}\label{se2eq1}
f_\bc(x) = \sum_{i=1}^{k_n} c_{1,1,i}^{(L)} \cdot f_{i,i}^{(L)}(x) + c_{1,1,0}^{(L)}
\end{equation}
for some $c_{1,1,0}^{(L)}, \dots, c_{1,1,k_n}^{(L)} \in \mathbb{R}$, where $f_{i,i}^{(L)}$ are recursively defined by
\begin{equation}
  \label{se2eq2}
f_{k,i}^{(r)}(x) = \sigma\left(\sum_{j=1}^{k_{0}} c_{k,i,j}^{(r-1)} \cdot f_{k,j}^{(r-1)}(x) + c_{k,i,0}^{(r-1)} \right)
\end{equation}
for some $c_{k,i,0}^{(r-1)}, \dots, c_{k,i, k_{0}}^{(r-1)} \in \mathbb{R}$
$(r=2, \dots, L)$
and
\begin{equation}
  \label{se2eq3}
f_{k,i}^{(1)}(x) = \sigma \left(\sum_{j=1}^d c_{k,i,j}^{(0)} \cdot x^{(j)} + c_{k,i,0}^{(0)} \right)
\end{equation}
for some $c_{k,i,0}^{(0)}, \dots, c_{k,i,d}^{(0)} \in \mathbb{R}$.

We learn the weight vector
$\bc=(c_{k,i,j}^{(s)})_{k,i,j,s}$
of our neural nework by the gradient descent. We initialize
$\bc^{(0)}$ by setting
\begin{equation}
  \label{se2eq5}
  c_{1,1,k}^{(L)}=0 \quad \mbox{for } k=0, \dots, k_n,
\end{equation}
and by choosing all others weights
randomly such that all weights $c_{k,i,j}^{(s)}$ with $s<L$
are independent
uniformly
distributed on $[-n^4, n^4]$, and we set
\[
\bc^{(t+1)}=\bc^{(t)} - \lambda_n \cdot (\nabla_\bc F_n)(\bc^{(t)})
\quad (t=0, \dots, t_n-1)
\]
where
\[
F_n(\bc) = \frac{1}{n} \sum_{i=1}^n |f_\bc(X_i)-Y_i|^2
\]
is the empirical $L_2$ risk of the network $f_\bc$ on the training data.
The the step size $\lambda_n>0$ and the number $t_n$
of gradient descent steps will be chosen below.

Because of (\ref{se2eq5}) we have
\[
F_n(\bc^{(0)}) = \frac{1}{n} \sum_{i=1}^n |Y_i|^2.
\]

\section{Main results}
\label{se3}
	Our first result shows that our estimate is able
        to achieve with high probability a very small error
        on the training data  in case  that $k_n$, $\lambda_n$ and
$t_n$ are suitably chosen.
	\begin{theorem}
	  \label{th1}
          Let $k_0 \in \N$ with $k_0 \geq 2 \cdot d$, let $L \in \N$ with
          $L \geq 2$, set
          \[
k_n=n^{5 \cdot (L-2) \cdot (k_0^2+k_0) + 5 \cdot k_0 \cdot (d+2)+7},
\]
\[
\lambda_n= \frac{1}{n^{8  (L-2) \cdot (k_0^2+k_0) + 8 \cdot k_0 \cdot (d+2) + 16L +15}}
\]
and
\[
t_n = 2 \cdot n^{
8 \cdot (L-2) \cdot (k_0^2+k_0) + 8 \cdot k_0 \cdot (d+2)+ 16L +17
},
\]
and define the estimate as in Section \ref{se2}.
Then for sufficiently large $n$ we have on the event
\[
\inf\{ \|X_i-X_j\|_\infty \, : \, 1 \leq i , j \leq n, X_i \neq X_j \}
\geq \frac{1}{(n+1)^3},
\]
\[
\max\{\|X_i\|_\infty \, : \, 1 \leq i \leq n \}\leq 1
\quad \mbox{and} \quad
\max\{|Y_i| \, : \, 1 \leq i \leq n \} \leq n^2
\]
 that
with probability at least $1-1/n$ the random choice of
$\bc^{(0)}$
leads to
\begin{equation}
\label{th1eq1}
\frac{1}{n} \sum_{i=1}^n | f_{\bc^{(t_n)}}(X_i)-Y_i|^2
\leq
\min_{g:\Rd \rightarrow \R} \frac{1}{n} \sum_{i=1}^n | g(X_i)-Y_i|^2
+
\frac{1}{n \cdot \log n}.
\end{equation}
	\end{theorem}

        \noindent
            {\bf Remark 1.} A corresponding result was shown
            in Kawaguchi and Huang (2019) for a fully connected
            network of
            much smaller size (linear instead of polynomial in
            the sample size as in Theorem \ref{th1} above),
            however there the learning rate
            of the gradient descent was set to zero for all
            weights $c_{k,i,j}^{(r)}$ with $r<L$. In contrast in
            our result the learning rate is positive for all
            weights.\\

            As our next result shows, any estimate which (as our estimate
            from Theorem \ref{th1}) achieves with high probability
            a very small error on the training data
            does in general not generalize well on a new independent data
            (provided we allow the distributions of $X$ which are concentrated
            on finite sets).

\begin{theorem}
		\label{th2}
Let $(X,Y)$, $(X_1,Y_1)$, \dots be independent
and identically distributed
$\Rd \times \R$-valued
random variables with $\EXP Y^2 < \infty$,
and let $U$ be an $\R^K$--valued  random variable independent of the
above
random variables.
Let $C_n$ be a subset of $\R^K$, and
let
\[
m_n(\cdot)=m_n(\cdot,(X_1,Y_1), \dots, (X_n,Y_n),U):\Rd \rightarrow \R
\]
be an estimate of $m$. Let $\kappa_n>0$ and let $\delta_n \leq
1/(n+1)^3$
and assume that $m_n$ satisfies
\[
\frac{1}{n} \sum_{i=1}^n | m_n(X_i)-Y_i|^2
\leq
\min_{g:\Rd \rightarrow \R}
\frac{1}{n} \sum_{i=1}^n | g(X_i)-Y_i|^2
+
\kappa_n
\]
whenever
\[
\inf \left\{
\|X_i-X_j\|_\infty \, : \,
1 \leq i,j \leq n, \, X_i \neq X_j
\right\}
\geq \delta_n
\quad \mbox{and}
\quad
U \in C_n.
\]
Then there exists a distribution of $(X,Y)$ such that $X \in [0,1]^d$ a.s.,
$Y \in \{-1,1\}$ a.s., $m(x)=0$ for all $x \in [0,1]^d$ and
such that we have for $n \geq 10$
\[
\EXP \int |m_n(x)-m(x)|^2 \PROB_X (dx) \geq
\frac{1}{5}
- n \cdot \kappa_n
-
\frac{1}{2} \cdot
\PROB_U(C_n^c).
\]
\end{theorem}

\begin{corollary}
\label{co1}
Let $p,C,c_1>0$ and let $\D^{(p,C)}$ be the class
of all distributions of $(X,Y)$ which satisfy
\begin{enumerate}
\item
$X \in [0,1]^d$ a.s.

\item $\sup_{x \in [0,1]^d} \EXP\{ Y^2 | X=x\} \leq c_1$

\item
$m(\cdot)=\EXP\{Y |X=\cdot\}$ is $(p,C)$--smooth.
\end{enumerate}
Let $m_n$ be the neural network regression estimate from Theorem
\ref{th1}.
Then we have for $n$ sufficiently large
\[
\sup_{(X,Y) \in \D^{(p,C)}} \EXP \int | m_n(x)-m(x)|^2 \PROB_X(dx)
\geq
\frac{1}{6}.
\]
\end{corollary}

\noindent
{\bf Proof.}
Let $U$ be the values for the random initialization of the weights
of the estimate in Theorem \ref{th1}.
By Theorem \ref{th1} we know that there exists a set $C_n$ of weights
such that (\ref{th1eq1}) holds for $n$ sufficiently large
whenever $U \in C_n$, where $\PROB_U(C^c_n) \leq 1/n$.
Hence the assumptions of Theorem \ref{th2} are satisfied with
$\kappa_n = 1/(n \cdot \log n)$.
Let $(X,Y)$ be the
distribution
from Theorem \ref{th2}. Then for $n$ sufficiently large
\begin{eqnarray*}
\sup_{(X,Y) \in \D^{(p,C)}} \EXP \int | m_n(x)-m(x)|^2 \PROB_X(dx)
&\geq&
\EXP \int | m_n(x)-m(x)|^2 \PROB_X(dx)
\\
&\geq& \frac{1}{5}
- n \cdot \frac{1}{n \cdot \log n}
-
\frac{1}{2} \cdot \frac{1}{n}
\\
&\geq&
\frac{1}{6}.
\end{eqnarray*}
\quad \hfill $\Box$

\noindent
    {\bf Remark 2.} Let $\D^{(p,C)}$ be the class of distributions
    of $(X,Y)$ introduced in Corollary \ref{co1}.
    It is well-known that there exist
    estimates $m_n$ which satisfy
    \[
    \sup_{(X,Y) \in \D^{(p,C)}} \EXP \int | m_n(x)-m(x)|^2 \PROB_X (dx)
    \leq
    c_2 \cdot n^{- \frac{2p}{2p+d}}
    \]
    and that no estimate can achieve a better rate of convergence
    (cf., Stone (1982) and Chapters 3 and 11 in Gy\"orfi et al. (2002)).
    Hence Corollary \ref{co1} implies that the estimate of Theorem \ref{th1}
    does not achieve the optimal minimax rate of convergence for the
    class $\D^{(p,C)}$, in fact the minimax $L_2$ error for this
    class does not even converge to zero, let alone in contrast to the optimal
    value.

\section{Proofs}
\label{se4}

\subsection{Proof of Theorem \ref{th1}}
\label{se4sub1}

\begin{lemma}
\label{le1}
Let $F:\R^K \rightarrow \R$ be differentiable, let $L_n>0$, set
\[
\lambda_n = \frac{1}{L_n},
\]
 let
$\ba_1\in \R^K$ and set
\[
\ba_2 =  \ba_1 - \lambda_n \cdot (\nabla_\ba F)(\ba_1).
\]
Then
\begin{equation}
\label{le1eq1}
\| (\nabla_\ba F)(\ba)-(\nabla_\ba F)(\ba_1) \|
\leq L_n \cdot \| \ba-\ba_1\|
\end{equation}
for all $\ba=\ba_1 + s \cdot (\ba_2-\ba_1)$, $s \in [0,1]$
implies
\[
F(\ba_2) \leq F(\ba_1) - \frac{1}{2 \cdot L_n} \cdot
\| (\nabla_\ba F)(\ba_1)\|^2.
\]
\end{lemma}

\noindent
{\bf Proof.} See proof of Lemma 1 in Braun, Kohler and Walk (2019).
\hfill $\Box$

Set
\[
F_n(\bc)= \frac{1}{n} \sum_{i=1}^n |f_{\bc}(X_i)-Y_i|^2
\]
where $f_{\bc}$ is defined by (\ref{se2eq1})--(\ref{se2eq3}).

\begin{lemma}
\label{le2}
Let $f_\bc$ be defined by (\ref{se2eq1})--(\ref{se2eq3}) and
assume that for any $i \in \{1, \dots, n\}$ there
exists $j_i \in \{1, \dots, k_n\}$ such that
\begin{equation}
\label{le2eq1}
f_{j_i,j_i}^{(L)}(X_i) \geq 1 - \frac{2}{n^2}
\quad \mbox{and} \quad
\sup_{t \in \{1, \dots,n\}, X_t \neq X_i}
f_{j_i,j_i}^{(L)}(X_t) \leq \frac{2}{n^2}
\end{equation}
hold. Then we have for any $n \geq 5$
\begin{eqnarray*}
  \| (\nabla_\bc F_n(\bc)) \|^2
  &\geq&
  \frac{1}{n} \cdot
\left(
  \frac{1}{n}
  \sum_{i=1}^n |f_{\bc}(X_i)-Y_i|^2
  -
  \min_{g:\Rd \rightarrow \R}
  \frac{1}{n}
  \sum_{i=1}^n |g(X_i)-Y_i|^2
  \right)
.
\end{eqnarray*}
\end{lemma}

\noindent
    {\bf Proof.} Set
\[
\bar{m}_n(x)
=
\frac{
\sum_{i=1}^n Y_i \cdot I_{\{X_i=x\}}
}{
\sum_{i=1}^n I_{\{X_i=x\}}
}
\quad (x \in \Rd),
\]
where we use the convention $0/0=0$. We have
\begin{eqnarray*}
&&
\frac{1}{n} \sum_{i=1}^n | f(X_i)-Y_i|^2
=
\frac{1}{n} \sum_{i=1}^n | f(X_i)-\bar{m}_n(X_i)|^2
+
\frac{1}{n} \sum_{i=1}^n | \bar{m}_n(X_i)-Y_i|^2,
\end{eqnarray*}
since
\begin{eqnarray*}
&&
\frac{1}{n} \sum_{i=1}^n (f(X_i)-\bar{m}_n(X_i))\cdot
(\bar{m}_n(X_i)-Y_i)
\\
&&
=
\frac{1}{n} \sum_{x \in \{X_1, \dots, X_n\}}
(f(x)-\bar{m}_n(x)) \cdot
\sum_{1 \leq i \leq n : X_i=x}
(\bar{m}_n(X_i)-Y_i)
=
0.
\end{eqnarray*}
This implies
\[
  \min_{g:\Rd \rightarrow \R}
  \frac{1}{n}
  \sum_{i=1}^n |g(X_i)-Y_i|^2
  =
\frac{1}{n} \sum_{i=1}^n | \bar{m}_n(X_i)-Y_i|^2  
\]
and
\[
\frac{1}{n} \sum_{i=1}^n | f(X_i)-\bar{m}_n(X_i)|^2
=
\frac{1}{n} \sum_{i=1}^n | f(X_i)-Y_i|^2
-
  \min_{g:\Rd \rightarrow \R}
  \frac{1}{n}
  \sum_{i=1}^n |g(X_i)-Y_i|^2
  \]
  for any $f: \Rd \rightarrow \R$.

   Next we observe
\begin{eqnarray*}
\| (\nabla_\bc F_n(\bc)) \|^2
&=&
\sum_{k,i,j,s} \left|
\frac{\partial}{\partial c_{k.j.i}^{(s)}} F_n(\bc)
\right|^2
\\
&\geq&
\sum_{i \in \{1, \dots, n\}, j_i \neq j_t \mbox{ for all } t<i} \left|
\frac{\partial}{\partial c_{1,1,j_i}^{(L)}} F_n(\bc)
\right|^2
\\
&=&
\sum_{i \in \{1, \dots, n\}, j_i \neq j_t \mbox{ for all } t<i} \left|
\frac{2}{n}
\cdot
\sum_{t=1}^n
(f_\bc(X_t)-Y_t)
\cdot
\frac{\partial}{\partial c_{1,1,j_i}^{(L)}} f_\bc(X_t)
\right|^2
\\
&=&
\sum_{i \in \{1, \dots, n\}, j_i \neq j_t \mbox{ for all } t<i}\left|
\frac{2}{n}
\cdot
\sum_{t=1}^n
(f_\bc(X_t)-Y_t)
\cdot
 f_{j_i,j_i}^{(L)}(X_t)
\right|^2
\\
&\geq&
\sum_{i \in \{1, \dots, n\}, j_i \neq j_t \mbox{ for all } t<i}
\Bigg(
\frac{1}{2} \cdot
 \left|
 \frac{2}{n}
\cdot
 \sum_{t \in \{1, \dots, n\}, X_t=X_i}
(f_\bc(X_i)-Y_t)
\cdot
 f_{j_i,j_i}^{(L)}(X_i)
\right|^2
\\
&&
\hspace*{2cm}
-
\left|
\frac{2}{n}
\cdot
\sum_{t \in \{1, \dots, n\}, X_t \neq X_i}
(f_\bc(X_t)-Y_t)
\cdot
 f_{j_i,j_i}^{(L)}(X_t)
\right|^2
\Bigg),
\end{eqnarray*}
where the last inequality followed from $b^2 \leq 2 (b-a)^2 + 2 a^2$
which implies
\[
a^2 \geq \frac{1}{2} b^2 - (b-a)^2 \quad (a,b \in \R).
\]
Using
\[
\sum_{t \in \{1, \dots, n\}, X_t=X_i}
(f_\bc (X_i)-Y_t)
=
| \{
1 \leq k \leq n \, : \, X_k=X_i
\}|
\cdot
(f_\bc (X_i)-\bar{m}_n(X_i)),
\]
\begin{eqnarray*}
  &&
\sum_{t \in \{1, \dots, n\}, X_t \neq X_i}
(f_\bc (X_t)-Y_t) \cdot
 f_{j_i,j_i}^{(L)}(X_t)
\\
&&
=
\sum_{t \in \{1, \dots, n\}, X_t \notin \{X_i,X_1, \dots, X_{t-1}\}}
| \{
1 \leq k \leq n \, : \, X_k=X_t
\}|
\cdot
(f_\bc (X_t)-\bar{m}_n(X_t)) \cdot
 f_{j_i,j_i}^{(L)}(X_t),
\end{eqnarray*}
(\ref{le2eq1}) and the inequality of Jensen we conclude
\begin{eqnarray*}
&&
\| (\nabla_\bc F_n(\bc)) \|^2
\\
&&
\geq
\frac{2}{n^2}
\cdot
\sum_{i \in \{1, \dots, n\}, j_i \neq j_t \mbox{ for all } t<i}
|\{1 \leq k \leq n \, : \, X_k=X_i\}|^2 \cdot
(f_\bc(X_i)-\bar{m}_n(X_i))^2
\cdot
\left(1 - \frac{2}{n^2} \right)^2
\\
&&
\hspace*{1cm}
-
4 \cdot
\sum_{i \in \{1, \dots, n\}, j_i \neq j_t \mbox{ for all } t<i}
\quad
\sum_{t \in \{1, \dots, n\}, X_t \notin \{X_i,X_1, \dots, X_{t-1}\}}
\\
&&
\hspace*{3cm}
\frac{
| \{
1 \leq k \leq n \, : \, X_k=X_t
\}|
}{n}
\cdot
|f_\bc(X_t)-\bar{m}_n(X_t)|^2
\cdot \frac{4}{n^{4}}
\\
&&
\geq
\frac{2}{n^2}
\cdot
\sum_{i \in \{1, \dots, n\}, j_i \neq j_t \mbox{ for all } t<i}
|\{1 \leq k \leq n \, : \, X_k=X_i\}| \cdot
(f_\bc(X_i)-\bar{m}_n(X_i))^2
\cdot
\left(1 - \frac{2}{n^2} \right)^2
\\
&&
\hspace*{1cm}
-
4 \cdot
\sum_{i \in \{1, \dots, n\}, j_i \neq j_t \mbox{ for all } t<i}
\quad
\sum_{t \in \{1, \dots, n\}, X_t \notin \{X_i,X_1, \dots, X_{t-1}\}}
\\
&&
\hspace*{3cm}
\frac{
| \{
1 \leq k \leq n \, : \, X_k=X_t
\}|
}{n}
\cdot
|f_\bc(X_t)-\bar{m}_n(X_t)|^2
\cdot \frac{4}{n^{4}}
\\
&&
\geq
\frac{2}{n}
\cdot
\left(1 - \frac{2}{n^2} \right)^2
\cdot
\frac{1}{n}
\sum_{t=1}^n
(f_\bc(X_t)-\bar{m}_n(X_t))^2
\\
&&
\hspace*{3cm}
-
4 \cdot n \cdot
\frac{4}{n^{4}}
\cdot
\frac{1}{n}
\cdot
\sum_{t=1}^n
|f_\bc(X_t)-\bar{m}_n(X_t)|^2
\\
&&
=
\left(
\frac{2}{n}
\cdot
\left(1 - \frac{2}{n^2} \right)^2
-
\frac{16}{n^{3 }}
\right)
\cdot
\frac{1}{n}
\cdot
\sum_{t=1}^n
|f_\bc(X_t)-\bar{m}_n(X_t)|^2
\\
&&
=
\left(
\frac{2}{n}
-
\frac{8}{n^{3}}
+
\frac{8}{n^{5}}
  -
  \frac{16}{n^{3}}
\right)
\cdot
\frac{1}{n}
\cdot
\sum_{t=1}^n
|f_\bc(X_t)-\bar{m}_n(X_t)|^2
\\
&&
\geq
\frac{1}{n} \cdot
\frac{1}{n}
 \cdot
\sum_{t=1}^n
(f_\bc(X_t)-\bar{m}_n(X_t))^2
\\
&&
=
\frac{1}{n} \cdot
\left(
  \frac{1}{n}
  \sum_{i=1}^n |f_{\bc}(X_i)-Y_i|^2
  -
  \min_{g:\Rd \rightarrow \R}
  \frac{1}{n}
  \sum_{i=1}^n |g(X_i)-Y_i|^2
  \right)
.
\end{eqnarray*}

\quad \hfill $\Box$

\begin{lemma}
\label{le3}
Define $\bc^{(t)}$ by
\[
\bc^{(t+1)}=\bc^{(t)} - \lambda_n \cdot (\nabla_\bc F_n)(\bc^{(t)})
\quad (t=0, \dots, t_n-1)
\]
for some fixed $\bc^{(0)}$ and
\[
\lambda_n = \frac{1}{L_n}.
\]
Assume that (\ref{le1eq1}) holds for $F=F_n$ and
all $\ba_1=\bc^{(t)}$ and $\ba_2=\bc^{(t+1)}$
and any $t \in \{0,1,\dots,t_n-1\}$.
Furthermore assume that (\ref{le2eq1}) holds for
all $\bc=\bc^{(t)}$ $(t \in \{0,1,\dots,t_n-1\})$.
Then we have for any $n \geq 5$
\begin{eqnarray*}
  &&
  F_n(\bc^{(t_n)})
  -
  \min_{g:\Rd \rightarrow \R}
  \frac{1}{n} \sum_{i=1}^n |g(X_i)-Y_i|^2
  \\
  &&\leq
  \left(
1-\frac{1}{2 \cdot n \cdot L_n}
\right)^{t_n}
\cdot
\left(
F_n(\bc^{(0)})
-
  \min_{g:\Rd \rightarrow \R}
  \frac{1}{n} \sum_{i=1}^n |g(X_i)-Y_i|^2
  \right)
.
\end{eqnarray*}
\end{lemma}

\noindent
{\bf Proof.}
Application of Lemma \ref{le1} and Lemma \ref{le2}
implies for any $t \in \{0,\dots,t_n-1\}$
\begin{eqnarray*}
 && F_n(\bc^{(t+1)})
  -
    \min_{g:\Rd \rightarrow \R}
    \frac{1}{n} \sum_{i=1}^n |g(X_i)-Y_i|^2
    \\
&&\leq
  F_n(\bc^{(t)})- \frac{1}{2 \cdot L_n}\cdot \| (\nabla_\bc F_n)(\bc^{(t)}) \|^2
  -
    \min_{g:\Rd \rightarrow \R}
  \frac{1}{n} \sum_{i=1}^n |g(X_i)-Y_i|^2
\\
&&\leq
\left(
1-\frac{1}{2 \cdot n \cdot L_n}
\right)
\cdot
\left(
F_n(\bc^{(t)}) -
    \min_{g:\Rd \rightarrow \R}
    \frac{1}{n} \sum_{i=1}^n |g(X_i)-Y_i|^2
    \right).
\end{eqnarray*}
From this we can conclude
\begin{eqnarray*}
&&  F_n(\bc^{(t_n)})
  -
    \min_{g:\Rd \rightarrow \R}
    \frac{1}{n} \sum_{i=1}^n |g(X_i)-Y_i|^2
    \\
&&\leq
\left(
1-\frac{1}{ 2 \cdot n \cdot L_n}
\right)
\cdot
\left(
F_n(\bc^{(t_n-1)})
-
    \min_{g:\Rd \rightarrow \R}
    \frac{1}{n} \sum_{i=1}^n |g(X_i)-Y_i|^2
    \right)
\\
&& \leq 
\left(
1-\frac{1}{ 2 \cdot n \cdot L_n}
\right)^2
\cdot
\left(
F_n(\bc^{(t_n-2)})
-
    \min_{g:\Rd \rightarrow \R}
    \frac{1}{n} \sum_{i=1}^n |g(X_i)-Y_i|^2
    \right)
\\
&&
\leq 
\dots
\\
&&
\leq 
\left(
1-\frac{1}{2 \cdot n \cdot L_n}
\right)^{t_n}
\cdot
\left(
F_n(\bc^{(0)})
-
    \min_{g:\Rd \rightarrow \R}
    \frac{1}{n} \sum_{i=1}^n |g(X_i)-Y_i|^2
    \right)
.
\end{eqnarray*}
\quad
\hfill $\Box$

\begin{lemma}
\label{le4}
Let $F:\R^K \rightarrow \R_+$ be differentiable, let $t_n \in \N$,
and let $L_n>0$
be such that
\begin{equation}
\label{le4eq1}
\| (\nabla_\ba F)(\ba) \|_\infty\leq L_n \cdot c_3 \cdot n^{c_4}
\quad \mbox{holds for all } \ba \mbox{ with } \|\ba\|_\infty
\leq 2 \cdot c_3 \cdot n^{c_4}
\end{equation}
and
\begin{equation}
\label{le4eq2}
\| (\nabla_\ba F)(\ba_1) -  (\nabla_\ba F)(\ba_2) \|
\leq L_n \cdot \|\ba_1-\ba_2\|
\end{equation}
holds for all  $\ba_1$, $\ba_2$ with $\|\ba_1\|_\infty \leq
3 \cdot c_3 \cdot n^{c_4}$ and
$\|\ba_2\|_\infty \leq
3 \cdot c_3 \cdot n^{c_4}$.
Let $\ba^{(0)}$ be such that
\begin{equation}
\label{le4eq3}
\|\ba^{(0)}\| \leq
c_3 \cdot n^{c_4}
\end{equation}
and
\begin{equation}
\label{le4eq4}
\sqrt{
\frac{2 \cdot t_n}{L_n} \cdot F( \ba^{(0)})
}
\leq
c_3 \cdot n^{c_4},
\end{equation}
and set
\[
\ba^{(t+1)}
=
\ba^{(t)}
-
\lambda_n \cdot
 (\nabla_\ba F)(\ba^{(t)})
\quad (t \in \{0,1,\dots, t_n-1\}),
\]
where
\[
\lambda_n=\frac{1}{L_n}.
\]
Then we have
\[
\| \ba^{(t)} \|_\infty \leq 2 \cdot c_3 \cdot n^{c_4}
\quad (t \in \{0,1,\dots,t_n\}.
\]
\end{lemma}

\noindent
{\bf Proof.} We show
\begin{equation}
\label{ple4eq1}
\| \ba^{(s)} \|_\infty \leq 2 \cdot c_3 \cdot n^{c_4}
\quad (s \in \{0,\dots,t\})
\end{equation}
for all $t \in \{0,1,\dots,t_n\}$ by induction.

For $t=0$ the assertion follows from (\ref{le4eq3}).
So assume that (\ref{ple4eq1}) holds for some
$t \in \{0,1,\dots,t_n-1\}$. Then this together with (\ref{le4eq1})
implies that we have
\[
\|\ba^{(t+1)}\|_\infty \leq \|\ba^{(t)}\|_\infty + \frac{1}{L_n} \cdot
\| (\nabla_\ba F)(\ba^{(t)}) \|_\infty
\leq 3 \cdot c_3 \cdot n^{c_4}.
\]
From this, the induction hypothesis and Lemma \ref{le1} we can
conclude
\[
0 \leq F(\ba^{(s)}) \leq F(\ba^{(s-1}) -
\frac{1}{2 \cdot L_n} \cdot
\| (\nabla_\ba F)(\ba^{(s-1)}) \|^2
\]
for all $s \in \{0,\dots,t+1\}$ which implies
\[
0 \leq F(\ba^{(t+1)}) \leq
F(\ba^{(0})-
\sum_{s=1}^{t+1}
\frac{1}{2 \cdot L_n} \cdot
\| (\nabla_\ba F)(\ba^{(s-1)}) \|^2.
\]
Consequently we have
\[
\sum_{s=1}^{t+1}
\frac{1}{2 \cdot L_n} \cdot
\| (\nabla_\ba F)(\ba^{(s-1)}) \|^2
\leq
F(\ba^{(0}),
\]
which implies
\begin{eqnarray*}
\|\ba^{(t+1)}\|_\infty
&\leq&
\|\ba^{(t+1)}\|\\
&\leq&
 \| \ba^{(0)}\| + \sum_{s=1}^{t+1} \frac{1}{L_n} \cdot
\| (\nabla_\ba F)(\ba^{(s-1)}) \|\\
&\leq&
 \| \ba^{(0)}\| +  \sqrt{ \frac{t+1}{ L_n}}
\cdot
\sqrt{
\sum_{s=1}^{t+1} \frac{1}{ L_n} \cdot
\| (\nabla_\ba F)(\ba^{(s-1)}) \|^2}
\\
&\leq&
\| \ba^{(0)}\| +  \sqrt{ \frac{t+1}{L_n} \cdot 2 \cdot F(\ba^{(0)})}
\\
& \leq & 2 \cdot c_3 \cdot n^{c_4},
\end{eqnarray*}
where the last inequality followed from (\ref{le4eq3}) and (\ref{le4eq4}).
\quad \hfill $\Box$

\begin{lemma}
\label{le5}
Let $\sigma$ be the logistic squasher.
Let $k_n \in \N$ and $k_0 \in \N$ with $2 \cdot  k_0 \geq d$.
Let
$\bc=(c_{k,i,j}^{(s)})_{k,i,j,s}$
and
$\bar{\bc}=(\bar{c}_{k,i,j}^{(s)})_{k,i,j,s}$
be weight vectors and define $f_{\bc}$ and $f_{\bar{\bc}}$ by
\begin{equation}\label{le5eq1}
f_\bc(x) = \sum_{i=1}^{k_n} c_{1,1,i}^{(L)} \cdot f_{i,i}^{(L)}(x) + c_{1,1,0}^{(L)}
\quad \mbox{and} \quad
f_{\bar\bc}(x) = \sum_{i=1}^{k_n} \bar{c}_{1,1,i}^{(L)} \cdot \bar{f}_{i,i}^{(L)}(x) + \bar{c}_{1,1,0}^{(L)}
\end{equation}
for $f_{i,i}^{(L)}$'s and $\bar{f}_{i,i}^{(L)}$'s recursively defined by
\begin{equation}
  \label{le5eq2}
f_{k,i}^{(r)}(x) = \sigma\left(\sum_{j=1}^{k_{0}} c_{k,i,j}^{(r-1)} \cdot f_{k,j}^{(r-1)}(x) + c_{k,i,0}^{(r-1)} \right)
\end{equation}
$(r=2, \dots, L)$
and
\begin{equation}
  \label{le5eq2b}
\bar{f}_{k,i}^{(r)}(x) = \sigma\left(\sum_{j=1}^{k_{0}} \bar{c}_{k,i,j}^{(r-1)} \cdot \bar{f}_{k,j}^{(r-1)}(x) + \bar{c}_{k,i,0}^{(r-1)} \right)
\end{equation}
$(r=2, \dots, L)$
and
\begin{equation}
  \label{le5eq3}
f_{k,i}^{(1)}(x) = \sigma \left(\sum_{j=1}^d c_{k,i,j}^{(0)} \cdot x^{(j)} + c_{k,i,0}^{(0)} \right)
\quad \mbox{and} \quad
\bar{f}_{k,i}^{(1)}(x) = \sigma \left(\sum_{j=1}^d
  \bar{c}_{k,i,j}^{(0)} \cdot x^{(j)}
+ \bar{c}_{k,i,0}^{(0)} \right).
\end{equation}
{\bf a)} For any $k \in \{1,\dots,k_n\}$ and any $x \in \Rd$ we have
\[
| f_{k,k}^{(L)}(x) - \bar{f}_{k,k}^{(L)}(x)| \leq  (2 \cdot
k_0+1)^L \cdot
(\max\{\|\bc\|_\infty, \|x\|_\infty,1\})^{L}
 \cdot
\max_{i,j,s:s<L}|c_{k,i,j}^{(s)}
-\bar{c}_{k,i,j}^{(s)}|.
\]
{\bf b)}
For any $x \in \Rd$ we have
\[
|f_{\bc}(x)-f_{\bar{\bc}}(x)|
\leq
(2 \cdot k_n+1)
\cdot
( 2 \cdot k_0+1)^L \cdot (\max\{\|\bc\|_\infty, \|x\|_\infty,1\})^{L+1}
 \cdot \|\bc-\bar{\bc}\|_{\infty}.
\]
\end{lemma}

\noindent
{\bf Proof.} {\bf a)} We show by induction
\begin{equation}
\label{ple5eq1}
| f_{r,k}^{(l)}(x) - \bar{f}_{r,k}^{(l)}(x)|
\leq
(2 \cdot k_0+1)^l \cdot (\max\{\|\bc\|_\infty, \|x\|_\infty,1\})^{l}
 \cdot \max_{i,j,s:s<L}|c_{k,i,j}^{(s)}-\bar{c}_{k,i,j}^{(s)}|
\end{equation}
$(l \in \{1, \dots, L\})$.
The logistic squasher satisfies $|\sigma^\prime(x)|
= |\sigma(x) \cdot (1-\sigma(x))| \leq 1$, hence
it is Lipschitz continuous with Lipschitz constant one.
This implies
\begin{eqnarray*}
\left|
f_{k,i}^{(1)}(x)
-
\bar{f}_{k,i}^{(1)}(x)
\right|
&
\leq
&
\sum_{j=1}^d |c_{k,i,j}^{(0)} - \bar{c}_{k,i,j}^{(0)}| \cdot |x^{(j)}|
+
|c_{k,i,0}^{(0)} - \bar{c}_{k,i,0}^{(0)}|
\\
&
\leq &
(2 \cdot k_0 +1) \cdot \max\{
\|x\|_\infty,1\} \cdot \max_{i,j,s:s<L}|c_{k,i,j}^{(s)}-\bar{c}_{k,i,j}^{(s)}|.
\end{eqnarray*}
Assume now that (\ref{ple5eq1}) holds for some $r-1$, where
$r \in \{2, \dots, L\}$. Then
\begin{eqnarray*}
&&
\left| f_{k,i}^{(r)}(x)
-
\bar{f}_{k,i}^{(r)}(x)
\right|
\\
&&
\leq
\sum_{j=1}^{k_{0}} |c_{k,i,j}^{(r-1)}|  \cdot
   |f_{k,j}^{(r-1)}(x)- \bar{f}_{k,j}^{(r-1)}(x)|
+
\sum_{j=1}^{k_{0}} |c_{k,i,j}^{(r-1)}- \bar{c}_{k,i,j}^{(r)}(x)|  \cdot
   |\bar{f}_{k,j}^{(r-1)}(x)|
 \\
&&
\quad
+
|c_{k,i,0}^{(r-1)} - \bar{c}_{k,i,0}^{(r-1)}|
\\
&&
\leq
k_0 \cdot \|\bc\|_\infty \cdot \max_{j=1, \dots, k_0}
|f_{k,j}^{(r-1)}(x)- \bar{f}_{k,j}^{(r-1)}(x)|
+
(k_0+1) \cdot \max_{i,j,s:s<L}|c_{k,i,j}^{(s)}-\bar{c}_{k,i,j}^{(s)}|
\\
&&
\leq
(2 k_0 +1) \cdot \max\{\|\bc\|_\infty,\|x\|_\infty,1\}
\\
&&
\hspace*{2cm}
\cdot
\max \left\{
\max_{j=1, \dots, k_0}
|f_{k,j}^{(r-1)}(x)- \bar{f}_{k,j}^{(r-1)}(x)|,
 \max_{i,j,s:s<L}|c_{k,i,j}^{(s)}-\bar{c}_{k,i,j}^{(s)}|
\right\}
\\
&&
\leq
  (2 \cdot k_0+1)^r \cdot (\max\{\|\bc\|_\infty, \|x\|_\infty,1\})^{r}
 \cdot \max_{i,j,s:s<L}|c_{k,i,j}^{(s)}-\bar{c}_{k,i,j}^{(s)}|.
\end{eqnarray*}

\noindent
{\bf b)}
Because of
\begin{eqnarray*}
&&
|f_{\bc}(x)-\bar{f}_{\bc}(x)|
\\
&&
=
\left|
\sum_{i=1}^{k_n} c_{1,1,i}^{(L)} \cdot f_{i,i}^{(L)}(x) + c_{1,1,0}^{(L)}
-
\sum_{i=1}^{k_n} \bar{c}_{1,1,i}^{(L)} \cdot \bar{f}_{i,i}^{(L)}(x) - \bar{c}_{1,1,0}^{(L)}
\right|
\\
&&
\leq
\left|
\sum_{i=1}^{k_n} c_{1,1,i}^{(L)}\cdot (f_{i,i}^{(L)}(x) - \bar{f}_{i,i}^{(L)}(x))
\right|
+
\left|
\sum_{i=1}^{k_n} (c_{1,1,i}^{(L)}-\bar{c}_{1,1,i}^{(L)}) \cdot
   \bar{f}_{i,i}^{(L)}(x)
+
 c_{1,1,0}^{(L)}-\bar{c}_{1,1,0}^{(L)}
\right|
\\
&&
\leq
k_n \cdot
\max_{i} |c_{1,1,i}^{(L)}| \cdot \max_{i} |f_{i,i}^{(L)}(x) -
   \bar{f}_{i,i}^{(L)}(x)|
+
(k_n+1) \cdot \max_{i}|c_{1,1,i}^{(L)}-\bar{c}_{1,1,i}^{(L)}|,
\end{eqnarray*}
the assertion follows from a).
\hfill $\Box$

\begin{lemma}
\label{le6}
Let $\sigma$ be the logistic squasher.
Define $f_{\bc}$ by (\ref{se2eq1})-(\ref{se2eq3}) and set
\[
F_n(\bc) = \frac{1}{n} \sum_{i=1}^n |f_\bc(X_i)-Y_i|^2.
\]
Let $c_3,c_4 \geq 1$.
Assume
$\|\bc_1\|_\infty \leq c_3 \cdot n^{c_4}$,
$\|\bc_2\|_\infty \leq c_3 \cdot n^{c_4}$
and
\[
\max_{i=1, \dots,n} \|X_i\|_\infty \leq c_3 \cdot n^{c_4}
\quad \mbox{and} \quad
\max_{i=1, \dots,n} |Y_i| \leq c_3 \cdot n^{c_4}.
\]
Set
\[
L_n
=
45 \cdot L \cdot 3 ^L
\cdot \left(
\max\{k_0,L,d\}
\right)^{3/2}
\cdot k_0^{2L}
\cdot k_n^{3/2}  \cdot
\left( c_3 \cdot n^{c_4} \right)^{4L+1}.
\]
Then we have
\begin{equation}
\label{le6eq1}
\| (\nabla_\bc F_n)(\bc_1) \|_\infty
\leq L_n \cdot c_3 \cdot n^{c_4}
\end{equation}
and
\begin{equation}
\label{le6eq2}
\| (\nabla_\bc F_n)(\bc_1) -  (\nabla_\bc F_n)(\bc_2) \|
\leq L_n \cdot \|\bc_1-\bc_2\|.
\end{equation}
\end{lemma}

\noindent
    {\bf Proof.}
    {\it In the first step of the proof} we compute the partial
    derivatives of $F_n(\bc)$.
    We have
    \[
    \frac{\partial}{\partial c_{k,i,j}^{(r)}} F_n(\bc)
    =
\frac{2}{n} \sum_{l=1}^n
\left( f_\bc(X_l)-Y_l \right)
\cdot
\frac{\partial f_\bc}{\partial c_{k,i,j}^{(r)}} (X_l).
    \]
The recursive definition of $f_\bc$ together with the chain rule imply
\begin{eqnarray*}
&&
\frac{\partial f_\bc}{\partial
   c_{1,1,i}^{(L)}} (X_l)
=
 f_{i,i}^{(L)}(X_l)
\end{eqnarray*}
(where we have set $f_{0,0}^{(L)}(x)=1$)
and in case $\bar{r}<L$
\[
\frac{\partial f_\bc}{\partial
   c_{\bar{k},\bar{i},\bar{j}}^{(\bar{r})}} (X_l)
=
\sum_{i=1}^{k_n} c_{1,1,i}^{(L)}\cdot
\frac{\partial f_{i,i}^{(L)}}{\partial
  c_{\bar{k},\bar{i},\bar{j}}^{(\bar{r})}} (X_l)
=
c_{1,1,\bar{k}}^{(L)}\cdot
\frac{\partial f_{\bar{k},\bar{k}}^{(L)}}{\partial
   c_{\bar{k},\bar{i},\bar{j}}^{(\bar{r})}} (X_l).
\]
In case $0 \leq \bar{r}<r$ and $r>1$ we have
\begin{eqnarray*}
  &&
\frac{\partial f_{k,i}^{(r)}}{\partial
   c_{k,\bar{i},\bar{j}}^{(\bar{r})}} (X_l)
\\
&=&
\sigma^\prime \left(\sum_{j=1}^{k_{0}} c_{k,i,j}^{(r-1)} \cdot
    f_{k,j}^{(r-1)}(X_l) + c_{k,i,0}^{(r-1)} \right)
\\
&&
\hspace*{3cm}
\cdot
\frac{\partial }{\partial
   c_{k,\bar{i},\bar{j}}^{(\bar{r})}}
\left(\sum_{j=1}^{k_{0}} c_{k,i,j}^{(r-1)} \cdot
f_{k,j}^{(r-1)}(X_l) + c_{k,i,0}^{(r-1)} \right)
\\
&=&
\sigma \left(\sum_{j=1}^{k_{0}} c_{k,i,j}^{(r-1)} \cdot
    f_{k,j}^{(r-1)}(X_l) + c_{k,i,0}^{(r-1)} \right)
\\
&&
\hspace*{3cm}
\cdot \left( 1 - \sigma \left(\sum_{j=1}^{k_{0}} c_{k,i,j}^{(r-1)} \cdot
    f_{k,j}^{(r-1)}(X_l) + c_{k,i,0}^{(r-1)} \right) \right)
\\
&&
\hspace*{3cm}
\cdot
\frac{\partial }{\partial
   c_{k,\bar{i},\bar{j}}^{(\bar{r})}}
\left(\sum_{j=1}^{k_{0}} c_{k,i,j}^{(r-1)} \cdot
f_{k,j}^{(r-1)}(X_l) + c_{k,i,0}^{(r-1)} \right).
\end{eqnarray*}
Next we explain how we can compute
\[
\frac{\partial }{\partial
   c_{k,\bar{i},\bar{j}}^{(\bar{r})}}
\left(\sum_{j=1}^{k_{0}} c_{k,i,j}^{(r-1)} \cdot
f_{k,j}^{(r-1)}(X_l) + c_{k,i,0}^{(r-1)} \right).
\]
In case $\bar{r}=r-1>0$ we have
\[
\frac{\partial }{\partial
   c_{k,\bar{i},\bar{j}}^{(r-1)}}
\left(\sum_{j=1}^{k_{0}} c_{k,i,j}^{(r-1)} \cdot
    f_{k,j}^{(r-1)}(X_l)+ c_{k,i,0}^{(r-1)} \right)
=
    f_{k,\bar{j}}^{(r-1)}(X_l) \cdot 1_{\{\bar{i}=i\}}
\]
(where we have set
$    f_{\bar{k},0}^{(r-1)}(x)=1$),
and in case $\bar{r} < r-1$ we get
\begin{eqnarray*}
&&
\frac{\partial }{\partial
   c_{k,\bar{i},\bar{j}}^{(\bar{r})}}
\left(\sum_{j=1}^{k_{0}} c_{k,i,j}^{(r-1)} \cdot
    f_{k,j}^{(r-1)}(X_l)+ c_{k,i,0}^{(r-1)} \right)
=
\sum_{j=1}^{k_{0}} c_{k,i,j}^{(r-1)} \cdot
\frac{\partial }{\partial
   c_{k,\bar{i},\bar{j}}^{(\bar{r})}}
    f_{k,j}^{(r-1)}(X_l).
\end{eqnarray*}

And in case
$r=2$ and $\bar{r}=0$ we have
\begin{eqnarray*}
  &&
\frac{\partial f_{k,i}^{(1)}}{\partial
   c_{k,\bar{i},\bar{j}}^{(0)}} (X_l)
\\
&&
=
\sigma^\prime \left(\sum_{j=1}^d c_{k,i,j}^{(0)} \cdot X_l^{(j)} + c_{k,i,0}^{(0)} \right)
\cdot
X_l^{(\bar{j})}
\cdot
1_{\{\bar{i}=i\}}
\\
&&
=
\sigma \left(\sum_{j=1}^d c_{k,i,j}^{(0)} \cdot X_l^{(j)} + c_{k,i,0}^{(0)} \right)
\cdot
\left(
1-
\sigma \left(\sum_{j=1}^d c_{k,i,j}^{(0)} \cdot X_l^{(j)} + c_{k,i,0}^{(0)} \right)
\right)
\cdot
X_l^{(\bar{j})}
\cdot
1_{\{\bar{i}=i\}},
\end{eqnarray*}
where we have set $X_l^{(0)}=1$.

{\it In the second step of the proof}
we show for $x \in \Rd$ with
$\|x\|_\infty \leq c_3 \cdot n^{c_4}$
and
$\bc$, $\bc_1$, $\bc_2$ with
$\|\bc\|_\infty \leq c_3 \cdot n^{c_4}$,
$\|\bc_1\|_\infty \leq c_3 \cdot n^{c_4}$
and
$\|\bc_2\|_\infty \leq c_3 \cdot n^{c_4}$,
\begin{equation}
  \label{ple6eqI}
  \left|
  \frac{\partial f_\bc (x)}{
\partial c_{k,i,j}^{(r)}
  }
  \right|
  \leq
  k_0^{L}
  \cdot
  \left( c_3 \cdot n^{c_4} \right)^{L+1}
  \end{equation}
and
\begin{equation}
  \label{ple6eqII}
  \left|
  \frac{\partial f_{\bc_1} (x)}{
\partial c_{k,i,j}^{(r)}
  }
  -
    \frac{\partial f_{\bc_2} (x)}{
\partial c_{k,i,j}^{(r)}
  }
\right|
  \leq
\bar{L}_n \cdot \| \bc_1-\bc_2\|_\infty,
    \end{equation}
where
\[
\bar{L}_n
=
4 L \cdot
3^L \cdot k_0^{2L-2} \cdot
 \left( c_3 \cdot n^{c_4} \right)^{4L}.
\]

It is easy to see that the first step of the proof implies
\begin{eqnarray}
 && \frac{\partial f_\bc (x)}{
\partial c_{k,i,j}^{(r)}
  }
  =
  \sum_{s_{r+1}=1}^{k_0} \dots \sum_{s_{L-2}=1}^{k_0}
  f_{k,j}^{(r)}(x)
  \cdot
  f_{k,i}^{(r+1)}(x) \cdot (1- f_{k,i}^{(r+1)}(x) )
  \nonumber \\
  && \quad
  \cdot
  c_{k,s_{r+1},i}^{(r+1)} \cdot
  f_{k,s_{r+1}}^{(r+2)}(x) \cdot (1- f_{k,s_{r+1}}^{(r+2)}(x) )
  \cdot
c_{k,s_{r+2},s_{r+1}}^{(r+2)} \cdot
  f_{k,s_{r+2}}^{(r+3)}(x) \cdot (1- f_{k,s_{r+2}}^{(r+3)}(x) )
   \nonumber \\
  && \quad
  \cdots
   c_{k,s_{L-2},s_{L-3}}^{(L-2)} \cdot
   f_{k,s_{L-2}}^{(L-1)}(x) \cdot (1- f_{k,s_{L-2}}^{(L-1)}(x) )
   \cdot
   c_{k,k,s_{L-2}}^{(L-1)} \cdot
   f_{k,k}^{(L)}(x) \cdot (1- f_{k,k}^{(L)}(x) )
   \nonumber \\
  && \quad
     \cdot
     c_{1,1,k}^{(L)},
   \label{ple6eq*}
  \end{eqnarray}
where we have used the abbreviations
\[
f_{k,j}^{(0)}(x)
=
\left\{
\begin{array}{ll}
  x^{(j)} & \mbox{if } j \in \{1,\dots,d\} \\
  1 & \mbox{if } j=0
\end{array}
\right.
\]
and
\[
f_{k,0}^{(r)}(x)=1.
\]
Because of
\[
f_{k,i}^{(r)}(x) \in [0,1] \quad \mbox{if } r>0
\]
and
\[
|f_{k,i}^{(0)}(x)| \leq c_3 \cdot n^{c_4}
\]
and
\[
\| \bc \|_\infty \leq c_3 \cdot n^{c_4}
\]
this implies (\ref{ple6eqI}).

Next we prove (\ref{ple6eqII}).
The right-hand side of (\ref{ple6eq*}) is a sum of at most
$k_0^{L-2}$ products, where each product contains at most
$3L+1$ factors. In the worst case from these
$3L+1$ factors $L$ are Lipschitz continuous functions with Lipschitz constant
bounded by one, which are bounded in absolute value
by $c_3 \cdot n^{c_4}$. And according to the proof of Lemma
\ref{le5} (cf., (\ref{ple5eq1})) the remaining $2L+1$ factors are
Lipschitz continuous functions with Lipschitz constant bounded
by
\[
(2k_0+1)^L \cdot ( c_3 \cdot n^{c_4})^L,
\]
which are bounded in absolute value by $c_3 \cdot n^{c_4}$.

If $g_1, \dots, g_s:\R \rightarrow \R$
are Lipschitz continuous functions with Lipschitz
constants $L_{g_1}$, \dots, $L_{g_s}$, then
\[
\prod_{l=1}^s g_l
\quad \mbox{and} \quad
\sum_{l=1}^s g_l
\]
are Lipschitz continuous functions with Lipschitz
constant bounded by
\[
\sum_{l=1}^s L_{g_l} \cdot \prod_{k \in \{1, \dots, s\} \setminus \{l\}}
\| g_k\|_\infty
\leq
s \cdot \max_l  L_{g_l}  \cdot (\max_k \| g_k\|_\infty)^{s-1}
\]
and by
\[
\sum_{l=1}^s L_{g_l} \leq s \cdot  \max_l  L_{g_l},
\]
respectively.
This implies that (\ref{ple6eq*}) is Lipschitz
continuous with Lipschitz constant bounded by
\[
k_0^{L-2} \cdot (3L+1) \cdot (2k_0+1)^L \cdot (c_3 n^{c_4})^L
\cdot (c_3 n^{c_4})^{3L}.
\]

In the {\it third step of the proof} we show (\ref{le6eq1}).
We have
\begin{eqnarray*}
\| (\nabla_\bc F_n)(\bc) \|_\infty
&=&
\max_{k,i,j,r}
\left|
\frac{2}{n} \sum_{l=1}^n
\left( f_\bc(X_l)-Y_l \right)
\cdot
\frac{\partial f_\bc}{\partial c_{k,i,j}^{(r)}} (X_i)
\right|
\\
&\leq&
2 \cdot \left(
(k_n+1) \cdot \|\bc\|_\infty + \max_{i=1,\dots,n}|Y_i|
\right)
\cdot
\max_{l,k,i,j,r}
\left|
\frac{\partial f_\bc}{\partial c_{k,i,j}^{(r)}} (X_l)
\right|
\\
&\leq&
6 \cdot k_n \cdot c_3 n^{c_4} \cdot
\max_{l,k,i,j,r}
\left|
\frac{\partial f_\bc}{\partial c_{k,i,j}^{(r)}} (X_l)
\right|
.
\end{eqnarray*}
From this the result follows by (\ref{ple6eqI}).

{\it In the fourth step of the proof} we show (\ref{le6eq2}).
Because of
\[
\| (\nabla_\bc F_n)(\bc_1) -  (\nabla_\bc F_n)(\bc_2) \|
=
\left(
\sum_{k,i,j,r}
\left|
\frac{\partial F_n}{\partial c_{k,i,j}^{(r)}} (\bc_1)
-
\frac{\partial F_n}{\partial c_{k,i,j}^{(r)}} (\bc_2)
\right|^2
\right)^{1/2}
\]
and
\[
\frac{\partial F_n}{\partial c_{k,i,j}^{(r)}} (\bc)
=
\frac{2}{n} \sum_{l=1}^n
\left( f_\bc(X_l)-Y_l \right)
\cdot
\frac{\partial f_\bc}{\partial c_{k,i,j}^{(r)}} (X_l)
\]
we have
\begin{eqnarray}
  &&
  \| (\nabla_\bc F_n)(\bc_1) -  (\nabla_\bc F_n)(\bc_2) \|
  \nonumber \\
  &&
\leq
\sqrt{
  k_n \cdot
  (k_0+1+(L-2) \cdot (k_0^2+k_0) + k_0 \cdot (d+1))+k_n+1
}
\nonumber \\
&&
\hspace*{1cm}
\cdot 2 \cdot
\max_{k,i,j,r,l}
\left|
\left( f_{\bc_1}(X_l)-Y_l \right)
\cdot
\frac{\partial f_{\bc_1}}{\partial c_{k,i,j}^{(r)}} (X_l)
-
\left( f_{\bc_2}(X_l)-Y_l \right)
\cdot
\frac{\partial f_{\bc_2}}{\partial c_{k,i,j}^{(r)}} (X_l)
\right|.
\label{ple6eq**}
\end{eqnarray}
By Lemma \ref{le5} we know
\begin{equation}
\label{ple6eqIII}
|f_{\bc_1}(X_l) - f_{\bc_2}(X_l)| \leq
(2 k_n+1) \cdot (2 k_0 +1)^L \cdot (c_3 n^{c_4})^{L+1}
\cdot \| \bc_1 - \bc_2 \|_\infty.
\end{equation}
Trivially,
\begin{equation}
\label{ple6eqIV}
|f_{\bc}(X_l)-Y_l| \leq (k_n+1) \cdot c_3 n^{c_4} +  c_3 n^{c_4}
    = (k_n+2) \cdot c_3 n^{c_4}.
\end{equation}
If
$g_i$ are Lipschitz continuous functions with Lipschitz constants
$L_{g_i}$,
then
 $g_1 \cdot g_2$
is Lipschitz continuous with Lipschitz constant
\[
\|g_1\|_\infty \cdot L_{g_2} + \|g_2\|_\infty \cdot L_{g_1}.
\]
Combining this with (\ref{ple6eqI}), (\ref{ple6eqII}), (\ref{ple6eqIII}) and (\ref{ple6eqIV})
we get that
\[
\bc
\mapsto
\left( f_\bc(X_l)-Y_l \right)
\cdot
\frac{\partial f_\bc}{\partial c_{k,i,j}^{(r)}} (X_l)
\]
is Lipschitz continuous with Lipschitz constant
bounded by
\begin{eqnarray*}
  &&
  (k_n+2) \cdot c_3 n^{c_4}
  \cdot 4L \cdot
3^L \cdot k_0^{2L-2} \cdot
\left( c_3 \cdot n^{c_4} \right)^{4L}
\\
&&
\hspace*{3cm}
+
k_0^L \cdot (c_3 n^{c_4})^{L+1}
\cdot
(2 k_n+1) \cdot (2 k_0 +1)^L \cdot (c_3 n^{c_4})^{L+1}
\\
&&
\leq
15 \cdot k_n \cdot L \cdot 3^L \cdot k_0^{2L}  \cdot (c_3 n^{c_4})^{4L+1}.
  \end{eqnarray*}
This together with (\ref{ple6eq**}) implies the assertion.
 \hfill $\Box$

 \begin{lemma}
   \label{le7}
   Let $\sigma$ be the logistic squasher and let $n,d,k_0,L \in \N$ with
   $k_0 \geq 2 \cdot d$ and $L \geq 2$.
   Define $f_{1,1}^{(L)}: \R \rightarrow \R$ recursively by
   \[
   f_{1,k}^{(r)}(x)
   =
   \sigma \left(
\sum_{j=1}^{k_0} c_{1,k,j}^{(r-1)} \cdot f_{1,j}^{(r-1)}(x) + c_{1,k,0}^{(r-1)}
   \right)
   \]
   for some $c_{1,k,0}^{(r-1)}$, \dots, $c_{1,k,k_0}^{(r-1)} \in \R$
   $(r=2, \dots, L)$ and
   \[
   f_{1,k}^{(1)}(x)= \sigma \left(
\sum_{j=1}^d c_{1,k,j}^{(0)} \cdot x^{(j)} + c_{1,k,0}^{(0)}
   \right)
   \]
   for some $c_{1,k,0}^{(0)}$, \dots, $c_{1,k,d}^{(0)} \in \R$.
   Let  $\delta>0$ and let $\ba$, $\bb \in \Rd$ such that
\[
b^{(l)}-a^{(l)} \geq 2 \cdot \delta \quad \mbox{for all } l \in \{1,\dots,d\}.
\]
Assume
\begin{equation}
  \label{le7eq1}
  c_{1,1,1}^{(L-1)} \leq -4 \cdot (n+1),
\end{equation}
\begin{equation}
  \label{le7eq2}
  | c_{1,1,j}^{(L-1)}-c_{1,1,1}^{(L-1)}| \leq \frac{1}{2 k_0} \quad \mbox{for } j=2, \dots,d,
\end{equation}
\begin{equation}
  \label{le7eq2b}
  | c_{1,1,j}^{(L-1)}| \leq \frac{1}{2 k_0} \quad \mbox{for } j=2d+1, \dots,k_0,
\end{equation}
\begin{equation}
  \label{le7eq3}
  |c_{1,k,0}^{(L-1)} + \frac{1}{2}  \cdot c_{1,1,1}^{(L-1)}| \leq \frac{1}{2}
\quad \mbox{for } k \in \{1, \dots, 2d\},
\end{equation}
\begin{equation}
  \label{le7eq4}
  c_{1,k,k}^{(r-1)} \geq 8 \cdot \log(8d-1)  \quad \mbox{for } k \in \{1, \dots, 2d\} \mbox{ and }
  r \in \{2, \dots, L-1\},
\end{equation}
\begin{equation}
  \label{le7eq5}
  | c_{1,k,0}^{(r-1)} + \frac{1}{2} \cdot c_{1,k,k}^{(r-1)}| \leq \frac{\log (8d-1)}{k_0} \quad \mbox{for } k \in \{1, \dots, 2d\} \mbox{ and }
  r \in \{2, \dots, L-1\},
\end{equation}
\begin{equation}
  \label{le7eq6}
  |c_{1,k,j}^{(r-1)}| \leq \frac{\log(8d-1)}{k_0} \quad \mbox{ for }
  j \in \{1, \dots, k_0\} \setminus \{k\}, k \in \{1, \dots, 2d\}, r \in \{2, \dots, L-1\},
\end{equation}
\begin{equation}
  \label{le7eq7}
  c_{1,k,k}^{(0)} \leq - \frac{2}{\delta} \cdot \log(8d-1)  \quad \mbox{ for } k \in \{1, \dots, d\},
\end{equation}
\begin{equation}
  \label{le7eq8}
  |c_{1,k,0}^{(0)} + a^{(k)} \cdot c_{1,k,k}^{(0)}| \leq \frac{\log (8d-1)}{d}
  \quad \mbox{ for } k \in \{1, \dots, d\},
\end{equation}
\begin{equation}
  \label{le7eq8b}
  |c_{1,k,j}^{(0)}| \leq  \frac{\log(8d-1)}{d}  \quad \mbox{ for } k \in \{1, \dots, d\}, j \in \{1, \dots, d\} \setminus \{k\}
\end{equation}
\begin{equation}
  \label{le7eq9}
  c_{1,d+k,k}^{(0)} \geq \frac{2}{\delta} \cdot \log(8d-1)  \quad \mbox{ for } k \in \{1, \dots, d\},
\end{equation}
\begin{equation}
  \label{le7eq10}
  |c_{1,d+k,0}^{(0)} + b^{(k)} \cdot c_{1,d+k,k}^{(0)}| \leq \frac{\log (8d-1)}{d}
    \quad \mbox{ for } k \in \{1, \dots, d\}
\end{equation}
and
\begin{equation}
  \label{le7eq10b}
  |c_{1,d+k,j}^{(0)}| \leq  \frac{\log(8d-1)}{d}  \quad \mbox{ for } k \in \{1, \dots, d\}, j \in \{1, \dots, d\} \setminus \{k\}.
\end{equation}
Then $f_{1,1}^{(L)}$ satisfies for any $x \in [-1,1]^d$
\begin{equation}
\label{le7eq12}
f_{1,1}^{(L)}(x) \geq 1 - e^{-n} \quad \mbox{if }
x \in [a^{(1)}+\delta,b^{(1)}-\delta] \times \dots \times
 [a^{(d)}+\delta,b^{(d)}-\delta]
\end{equation}
and
\begin{equation}
\label{le7eq13}
f_{1,1}^{(L)}(x) \leq e^{-n} \quad \mbox{if }
x \notin [a^{(1)}-\delta,b^{(1)}+\delta] \times \dots \times
 [a^{(d)}-\delta,b^{(d)}+\delta].
\end{equation}
 \end{lemma}

 \noindent
     {\bf Proof.} Let $
x \in [a^{(1)}+\delta,b^{(1)}-\delta] \times \dots \times
[a^{(d)}+\delta,b^{(d)}-\delta]
\cap [-1,1]^d
 $.
 Then  we get for any $k \in \{1, \dots, d\}$
 by (\ref{le7eq7}), (\ref{le7eq8}) and (\ref{le7eq8b})
 \begin{eqnarray*}
   &&
   \sum_{j=1}^d c_{1,k,j}^{(0)} \cdot x^{(j)} + c_{1,k,0}^{(0)}
   \\
   &&
   =
   c_{1,k,k}^{(0)} \cdot ( x^{(k)} - a^{(k)})
   +
   c_{1,k,0}^{(0)} +  c_{1,k,k}^{(0)} \cdot a^{(k)}
   +
   \sum_{j \in \{1, \dots, d\} \setminus \{k\}} c_{1,k,j}^{(0)} \cdot x^{(j)}
   \\
   &&
   \leq
   - 2 \cdot \log(8d-1) +
   |c_{1,k,0}^{(0)} +  c_{1,k,k}^{(0)} \cdot a^{(k)}|
   +
   \sum_{j \in \{1, \dots, d\} \setminus \{k\}} |c_{1,k,j}^{(0)}|
   \\
   &&
   \leq
   - \log(8d-1).
 \end{eqnarray*}
 And by (\ref{le7eq9}), (\ref{le7eq10}) and (\ref{le7eq10b})
 we get for any $k \in \{1, \dots, d\}$
 \begin{eqnarray*}
   &&
   \sum_{j=1}^d c_{1,d+k,j}^{(0)} \cdot x^{(j)} + c_{1,d+k,0}^{(0)}
   \\
   &&
   =
   - c_{1,d+k,k}^{(0)} \cdot ( b^{(k)} - x^{(k)})
   +
   c_{1,d+k,0}^{(0)} +  c_{1,d+k,k}^{(0)} \cdot b^{(k)}
   +
   \sum_{j \in \{1, \dots, d\} \setminus \{k\}} c_{1,d+k,j}^{(0)} \cdot x^{(j)}
   \\
   &&
   \leq
   - 2 \cdot \log(8d-1) +
   |c_{1,d+k,0}^{(0)} +  c_{1,d+k,k}^{(0)} \cdot b^{(k)}|
   +
   \sum_{j \in \{1, \dots, d\} \setminus \{k\}} |c_{1,d+k,j}^{(0)}|
   \\
   &&
   \leq
   - \log(8d-1).
 \end{eqnarray*}
It is easy to see that
the logistic squasher satisfies
\begin{equation}
\label{ple7eq1}
\sigma(x) \geq 1 - \kappa
\quad \mbox{if} \quad
x \geq \log \left( \frac{1}{\kappa}-1 \right)
\quad \mbox{and} \quad
\sigma(x) \leq \kappa
\quad \mbox{if} \quad
x \leq - \log \left( \frac{1}{\kappa}-1 \right).
\end{equation}
Using this we get for any $k \in \{1, \dots, 2 d\}$
\[
f_{1,k}^{(1)}(x) \leq \sigma(-\log(8d-1))
=
\sigma \left(
- \log \left( \frac{1}{1/(8d)}-1 \right)
\right)
\leq \frac{1}{8d}
\leq \frac{1}{4}.
\]
Using (\ref{le7eq4}),  (\ref{le7eq5}) and  (\ref{le7eq6}), we can recursively conclude for $r=2, \dots, L-1$
that we have for any $k \in \{1, \dots, 2 d\}$
\begin{eqnarray*}
  &&
\sum_{j=1}^{k_0} c_{1,k,j}^{(r-1)} \cdot f_{1,j}^{(r-1)}(x) + c_{1,k,0}^{(r-1)}
\\
&&
=
c_{1,k,k}^{(r-1)} \cdot \left( f_{1,k}^{(r-1)}(x) - \frac{1}{2} \right)
+ c_{1,k,k}^{(r-1)} \cdot \frac{1}{2} + c_{1,k,0}^{(r-1)}
+
\sum_{j \in \{1, \dots, k_0\} \setminus \{k\}} c_{1,k,j}^{(r-1)} \cdot f_{1,j}^{(r-1)}(x)
\\
&&
\leq
- 2 \cdot \log(8d-1)
+ |c_{1,k,k}^{(r-1)} \cdot \frac{1}{2} + c_{1,k,0}^{(r-1)}|
+
\sum_{j \in \{1, \dots, k_0\} \setminus \{k\}} |c_{1,k,j}^{(r-1)}|
\\
&&
\leq - \log(8d-1)
\end{eqnarray*}
and
\[
f_{1,k}^{(r)}(x) \leq \sigma(-\log(8d-1))
\leq
 \frac{1}{8d}
\leq
 \frac{1}{4}.
 \]
 From this together with (\ref{le7eq1}), (\ref{le7eq2}), (\ref{le7eq2b})
 and (\ref{le7eq3}) we conclude
 \begin{eqnarray*}
   &&
   \sum_{j=1}^{k_0} c_{1,1,j}^{(L-1)} \cdot f_{1,j}^{(L-1)}(x) + c_{1,1,0}^{(L-1)}
   \\
   &&
   =
   c_{1,1,1}^{(L-1)} \cdot (\sum_{j=1}^{2d} f_{1,j}^{(L-1)}(x) -\frac{1}{2})
   + c_{1,1,0}^{(L-1)} +\frac{1}{2}  \cdot c_{1,1,1}^{(L-1)}
   \\
   &&
   \quad
   +
   \sum_{j=1}^{2d} (c_{1,1,j}^{(L-1)}-   c_{1,1,1}^{(L-1)}) \cdot f_{1,j}^{(L-1)}(x)
   +
   \sum_{j=2d+1}^{k_0} c_{1,1,j}^{(L-1)} \cdot f_{1,j}^{(L-1)}(x)
   \\
   &&
   \geq
   c_{1,1,1}^{(L-1)} \cdot (\sum_{j=1}^{2d} f_{1,j}^{(L-1)}(x) -\frac{1}{2})
   - |c_{1,k,0}^{(L-1)} +\frac{1}{2} \cdot c_{1,1,1}^{(L-1)}|
   \\
   &&
   \quad
   -
   \sum_{j=1}^{2d} |c_{1,1,j}^{(L-1)}-   c_{1,1,1}^{(L-1)}|
   -
   \sum_{j=2d+1}^{k_0} |c_{1,1,j}^{(L-1)}|
   \\
   &&
   \geq -4 \cdot
   (n+1) \cdot (2d \cdot \frac{1}{8d} - \frac{1}{2}) - \frac{1}{2}
   -  \sum_{j=1}^{2d} \frac{1}{2 k_0}
   -
   \sum_{j=2d+1}^{k_0} \frac{1}{2 k_0}
   \\
   &&
   \geq n \geq  \log(1/e^{-n} -1),
   \end{eqnarray*}
 which implies (\ref{le7eq12}).

In order to prove (\ref{le7eq13}) we assume that $x \in [-1,1]^d$
satisfies $ x^{(k)} \notin [a^{(k)}-\delta,b^{(k)}+\delta]$ for some
$k \in \{1, \dots, d\}$. In case $x^{(k)} < a^{(k)}-\delta$
we can argue similarly as above and conclude recursively
from (\ref{ple7eq1}) and (\ref{le7eq1})-(\ref{le7eq10b})
 \begin{eqnarray*}
   &&
   \sum_{j=1}^d c_{1,k,j}^{(0)} \cdot x^{(j)} + c_{1,k,0}^{(0)}
   \\
   &&
   =
   c_{1,k,k}^{(0)} \cdot ( x^{(k)} - a^{(k)})
   +
   c_{1,k,0}^{(0)} +  c_{1,k,k}^{(0)} \cdot a^{(k)}
   +
   \sum_{j \in \{1, \dots, d\} \setminus \{k\}} c_{1,k,j}^{(0)} \cdot x^{(j)}
   \\
   &&
   \geq
   2 \cdot \log(8d-1) -
   |c_{1,k,0}^{(0)} +  c_{1,k,k}^{(0)} \cdot a^{(k)}|
   -
   \sum_{j \in \{1, \dots, d\} \setminus \{k\}} |c_{1,k,j}^{(0)}|
   \\
   &&
   \geq
   \log(8d-1),
 \end{eqnarray*}
 which implies
 \[
 f_{1,k}^{(1)}(x) \geq \sigma( \log(8d-1)) = \sigma(\log(1/(1/(8d))-1))
 \geq
 1 - \frac{1}{8d}
 \geq \frac{3}{4}.
 \]
 Recursively we can conclude for $r=2, \dots, L-1$
\begin{eqnarray*}
  &&
\sum_{j=1}^{k_0} c_{1,k,j}^{(r-1)} \cdot f_{1,j}^{(r-1)}(x) + c_{1,k,0}^{(r-1)}
\\
&&
=
c_{1,k,k}^{(r-1)} \cdot \left( f_{1,k}^{(r-1)}(x) - \frac{1}{2} \right)
+ c_{1,k,k}^{(r-1)} \cdot \frac{1}{2} + c_{1,k,0}^{(r-1)}
+
\sum_{j \in \{1, \dots, k_0\} \setminus \{k\}} c_{1,k,j}^{(r-1)} \cdot f_{1,j}^{(r-1)}(x)
\\
&&
\geq
 2 \cdot \log(8d-1)
- |c_{1,k,k}^{(r-1)} \cdot \frac{1}{2} + c_{1,k,0}^{(r-1)}|
-
\sum_{j \in \{1, \dots, k_0\} \setminus \{k\}} |c_{1,k,j}^{(r-1)}|
\\
&&
\geq \log(8d-1)
\end{eqnarray*}
and
\[
f_{1,k}^{(r)}(x) \geq \sigma(\log(8d-1))
\geq
1- \frac{1}{8d}
\geq
 \frac{3}{4}.
 \]
 This yields
\begin{eqnarray*}
   &&
   \sum_{j=1}^{k_0} c_{1,1,j}^{(L-1)} \cdot f_{1,j}^{(L-1)}(x) + c_{1,1,0}^{(L-1)}
   \\
   &&
   =
   c_{1,1,1}^{(L-1)} \cdot (\sum_{j=1}^{2d} f_{1,j}^{(L-1)}(x) -\frac{1}{2})
   + c_{1,1,0}^{(L-1)} +\frac{1}{2}  \cdot c_{1,1,1}^{(L-1)}
   \\
   &&
   \quad
   +
   \sum_{j=1}^{2d} (c_{1,1,j}^{(L-1)}-   c_{1,1,1}^{(L-1)}) \cdot f_{1,j}^{(L-1)}(x)
   +
   \sum_{j=2d+1}^{k_0} c_{1,1,j}^{(L-1)} \cdot f_{1,j}^{(L-1)}(x)
   \\
   &&
   \leq
   c_{1,1,1}^{(L-1)} \cdot (\sum_{j=1}^{2d} f_{1,j}^{(L-1)}(x) -\frac{1}{2})
   + |c_{1,k,0}^{(L-1)} +\frac{1}{2}  \cdot c_{1,1,1}^{(L-1)}|
   \\
   &&
   \quad
   +
    \sum_{j=1}^{2d} |c_{1,1,j}^{(L-1)}-   c_{1,1,1}^{(L-1)}|
   +
   \sum_{j=2d+1}^{k_0} |c_{1,1,j}^{(L-1)}|
   \\
   &&
   \leq -4 \cdot
   (n+1) ( \frac{3}{4} - \frac{1}{2})  + \frac{1}{2}
   +  \sum_{j=1}^{2d} \frac{1}{2 k_0}
   +
   \sum_{j=2d+1}^{k_0} \frac{1}{2 k_0}
   \\
   &&
   \leq -n \leq  -\log(1/e^{-n} -1),
   \end{eqnarray*}
 which implies (\ref{le7eq13}).

In the same way we get the assertion in case
$x^{(k)} > b^{(k)}+\delta$. \hfill $\Box$

\noindent
    {\bf Remark 3.} It is easy to see that the number
    of weights of the neural network $f_{1,1}^{(L)}$
    is given by
    \[
(L-2) \cdot (k_0^2+k_0) + k_0 \cdot (d+2) +1.
    \]

\begin{lemma}
  \label{le8}
  Let $\sigma$ be the logistic squasher.
  Let $f_\bc$ be defined by (\ref{se2eq1})--(\ref{se2eq3}),
  let $k \in \{1, \dots, k_n\}$
  and assume that
\begin{equation}
  \label{le8eq1}
  \max_{t \in \{1, \dots, n\}}
  f_{k,k}^{(L)}(X_t) \cdot
(1-f_{k,k}^{(L)}(X_t))
  \leq e^{-n}
\end{equation}
holds.
Assume furthermore
$\|\bc\|_\infty \leq c_3 \cdot n^{c_4}$
and
\begin{equation}
\label{le8eq2}
\max_{j=1, \dots,n} \|X_j\|_\infty \leq c_3 \cdot n^{c_4}
\quad \mbox{and} \quad
\max_{j=1, \dots,n} |Y_j| \leq c_3 \cdot n^{c_4}.
\end{equation}
Set
\[
F_n(\bc) = \frac{1}{n} \sum_{i=1}^n |f_\bc(X_i)-Y_i|^2.
\]
Then we have for all $r<L$ and all $i$, $j$
\[
\left|
\frac{\partial F_n}{\partial c_{k,i,j}^{(r)}} (\bc)
\right|
\leq
2 \cdot \sqrt{F_n(\bc)} \cdot k_0^{L} \cdot (c_3 \cdot n^{c_4})^{L+1} \cdot e^{-n}
\]
\end{lemma}

\noindent
    {\bf Proof.} By the Cauchy-Schwarz inequality we get
    \begin{eqnarray*}
      \left|  \frac{\partial}{\partial c_{k,i,j}^{(r)}} F_n(\bc)
      \right|
      &=&
      \left|
\frac{2}{n} \sum_{l=1}^n
\left( f_\bc(X_l)-Y_l \right)
\cdot
\frac{\partial f_\bc}{\partial c_{k,i,j}^{(r)}} (X_l)
\right|
\\
&\leq&
2 \cdot \sqrt{F_n(\bc)} \cdot
\max_{l=1,\dots,n}
\left|
\frac{\partial f_\bc}{\partial c_{k,i,j}^{(r)}} (X_l)
\right|.
\end{eqnarray*}
    Using the recursive definition of $f_\bc$ together
    with
    (\ref{le8eq1}), $r<L$ and $\sigma^\prime(x)=\sigma(x) \cdot (1-\sigma(x))$
    we get
    \begin{eqnarray*}
      &&
      \left|
\frac{\partial f_\bc}{\partial c_{k,i,j}^{(r)}} (X_l)
\right|
\\
&&
=
      \left|
      \sum_{\bar{i}=1}^{k_n} c_{1,1,\bar{i}}^{(L)} \cdot
\frac{\partial f_{\bar{i},\bar{i}}^{(L)}}{\partial c_{k,i,j}^{(r)}}
      (X_l)
\right|
\\
&&
=
|c_{1,1,k}^{(L)}| \cdot
      \left|
\frac{\partial f_{k,k}^{(L)}}{\partial c_{k,i,j}^{(r)}}
      (X_l)
\right|
\\
&&
=
|c_{1,1,k}^{(L)}| \cdot
f_{k,k}^{(L)}(X_l)
\cdot
(1-f_{k,k}^{(L)}(X_l))
\cdot
\left|
\frac{\partial }{\partial
   c_{k,i,j}^{(r)}}
\left(\sum_{\bar{j}=1}^{k_{0}} c_{k,i,\bar{j}}^{(L-1)} \cdot
f_{k,\bar{j}}^{(L-1)}(X_l) + c_{k,i,0}^{(L-1)} \right)
\right|
\\
&&
\leq
|c_{1,1,k}^{(L)}| \cdot e^{-n} \cdot
\left|
\frac{\partial }{\partial
   c_{k,i,j}^{(r)} }
\left(\sum_{\bar{j}=1}^{k_{0}} c_{k,i,\bar{j}}^{(L-1)} \cdot
f_{k,\bar{j}}^{(L-1)}(X_l) + c_{k,i,0}^{(L-1)} \right)
\right|.
    \end{eqnarray*}
    As in the proof of Lemma \ref{le6} (cf., proof of (\ref{ple6eqI})) it is possible to show
    \[
    |c_{1,1,k}^{(L)}| \cdot
\left|
\frac{\partial }{\partial
   c_{k,i,j}^{(r)} }
\left(\sum_{\bar{j}=1}^{k_{0}} c_{k,i,\bar{j}}^{(L-1)} \cdot
f_{k,\bar{j}}^{(L-1)}(X_l) + c_{k,i,0}^{(L-1)} \right)
\right|
\leq k_0^{L} \cdot (c_3 \cdot n^{c_4})^{L+1},
\]
which implies the assertion.
    \hfill $\Box$

\noindent
    {\bf Proof of Theorem \ref{th1}.}
    The proof is divided into six steps.

    In the {\it first step of the proof} we show that
    for every $l \in \{1, \dots, n\}$ there exist (random)
    \[
    (\bar{c}_{1,i,j}^{(r)})_{i,j,r \, : r<L}
    \in \left[-n^4,
      n^4
    \right]^{
(L-2) \cdot (k_0^2+k_0) + k_0 \cdot (d+2) +1
      }
    \]
    such that for any
    $(c_{1,i,j}^{(r)})_{i,j,r \, : r<L}$ with
        \begin{equation}
      \label{pth1eq1}
      \max_{i,j,r:r<L} |c_{1,i,j}^{(r)}-\bar{c}_{1,i,j}^{(r)}| <
      \min \left\{
      2 (n+1), \frac{1}{16 k_0},\frac{1}{16},
      \frac{\log(8d-1)}{24 k_0}
      \right\}
        \end{equation}
    we have that any function $f_{1,1}^{(L)}$ corresponding
    to any $(\tilde{c}_{1,i,j}^{(r)})_{i,j,r \, : r<L}$ with
        \begin{equation}
      \label{pth1eq2}
      \max_{i,j,r:r<L} |\tilde{c}_{1,i,j}^{(r)}-c_{1,i,j}^{(r)}| <
\min \left\{
      2 (n+1), \frac{1}{16 k_0},\frac{1}{16},
      \frac{\log(8d-1)}{24 k_0}
      \right\}
       \end{equation}
    satisfies in case $\min\{ \|X_i-X_j\|_\infty : 1 \leq i,j \leq n, X_i \neq X_j\}\geq 1/(n+1)^3$ 
    \begin{equation}
      \label{pth1eq3}
      f_{1,1}^{(L)}(X_l) \geq 1 - e^{-n}
      \quad \mbox{and} \quad
    \max_{t \in \{1, \dots,n\}, \, X_t \neq X_l} f_{1,1}^{(L)}(X_t)  \leq e^{-n} .
    \end{equation}
    Set $\delta_n=1/(n+1)^3$ and
    $a^{(i)}=X_l^{(i)} - \frac{\delta_n}{2}$ and
    $b^{(i)}=X_l^{(i)} + \frac{\delta_n}{2}$ $(i=1, \dots,d)$.
    Then we have
    \[
    X_l \in \left[a^{(1)} + \frac{\delta_n}{4}, b^{(1)} - \frac{\delta_n}{4}
      \right]
    \times \dots \times
    \left[a^{(d)} + \frac{\delta_n}{4}, b^{(d)} - \frac{\delta_n}{4}
      \right],
    \]
    and
    \[
    \min\{ \|X_i-X_j\|_\infty : 1 \leq i,j \leq n, X_i \neq X_j\}\geq 1/(n+1)^3
    \]
    implies that we  also have
    \[
    X_t \notin \left[a^{(1)} - \frac{\delta_n}{4}, b^{(1)} + \frac{\delta_n}{4}
      \right]
    \times \dots \times
    \left[a^{(d)} - \frac{\delta_n}{4}, b^{(d)} + \frac{\delta_n}{4}
      \right]
    \]
    for all $t \in \{1, \dots,n\}$ with $X_t \neq X_l$.
    If
    $(\bar{c}_{1,i,j}^{(r)})_{i,j,r \, : r<L}$
    satisfies
\[
  \bar{c}_{1,1,1}^{(L-1)} \leq -8 \cdot (n+1),
\]
\[
  | \bar{c}_{1,1,j}^{(L-1)}-\bar{c}_{1,1,1}^{(L-1)}| \leq \frac{1}{4 k_0} \quad \mbox{for } j=2, \dots,d,
\]
\[
  | \bar{c}_{1,1,j}^{(L-1)}| \leq \frac{1}{4 k_0} \quad \mbox{for } j=2d+1, \dots,k_0,
\]
\[
  |\bar{c}_{1,k,0}^{(L-1)} + \frac{1}{2}  \cdot \bar{c}_{1,1,1}^{(L-1)}| \leq \frac{1}{4}
\quad \mbox{for } k \in \{1, \dots, 2d\},
\]
\[
  \bar{c}_{1,k,k}^{(r-1)} \geq 16 \cdot \log(8d-1)  \quad \mbox{for } k \in \{1, \dots, 2d\} \mbox{ and }
  r \in \{2, \dots, L-1\},
\]
\[
  | \bar{c}_{1,k,0}^{(r-1)} + \frac{1}{2} \cdot \bar{c}_{1,k,k}^{(r-1)}| \leq \frac{\log (8d-1)}{2 k_0} \quad \mbox{for } k \in \{1, \dots, 2d\} \mbox{ and }
  r \in \{2, \dots, L-1\},
\]
\[
  |\bar{c}_{1,k,j}^{(r-1)}| \leq \frac{\log(8d-1)}{2 k_0} \quad \mbox{ for }
  j \in \{1, \dots, k_0\} \setminus \{k\}, k \in \{1, \dots, 2d\}, r \in \{2, \dots, L-1\},
\]
\[
  \bar{c}_{1,k,k}^{(0)} \leq - \frac{4}{\delta_n} \cdot \log(8d-1)  \quad \mbox{ for } k \in \{1, \dots, d\},
\]
\[
  |\bar{c}_{1,k,0}^{(0)} + a^{(k)} \cdot \bar{c}_{1,k,k}^{(0)}| \leq \frac{\log (8d-1)}{2d}
  \quad \mbox{ for } k \in \{1, \dots, d\},
\]
\[
  |\bar{c}_{1,k,j}^{(0)}| \leq  \frac{\log(8d-1)}{2d}  \quad \mbox{ for } k \in \{1, \dots, d\}, j \in \{1, \dots, d\} \setminus \{k\}
\]
\[
  \bar{c}_{1,d+k,k}^{(0)} \geq \frac{4}{\delta_n} \cdot \log(8d-1)  \quad \mbox{ for } k \in \{1, \dots, d\},
\]
\[
|\bar{c}_{1,d+k,0}^{(0)} + b^{(k)} \cdot \bar{c}_{1,d+k,k}^{(0)}| \leq \frac{\log (8d-1)}{2 d}
\quad \mbox{for } k \in \{1,\dots,d\}
\]
and
\[
  |\bar{c}_{1,d+k,j}^{(0)}| \leq  \frac{\log(8d-1)}{2d}  \quad \mbox{ for } k \in \{1, \dots, d\}, j \in \{1, \dots, d\} \setminus \{k\},
\]
then it is easy to see that for any
$(c_{1,i,j}^{(r)})_{i,j,r \, : r<L}$ which satisfies
(\ref{pth1eq1}) we have that any
$(\tilde{c}_{1,i,j}^{(r)})_{i,j,r \, : r<L}$
which
satisfies (\ref{pth1eq2}) also
satisfies (\ref{le7eq1})-(\ref{le7eq10b}). Application of Lemma \ref{le7}
yields (\ref{pth1eq3}).

In the {\it second step of the proof} we show that
for $n$ sufficiently large
with
probability at least $1-n \cdot e^{-n}$ the weights in the random
initialization of the weights are chosen such that for
each $l \in \{1, \dots, n\}$ the weights for some index $k_l$
satisfy (\ref{pth1eq1}) (and hence all functions with weights
satisfying (\ref{pth1eq2}) satisfy (\ref{pth1eq3})).
We assume in the sequel that $n$ is sufficiently large.
If we sample the weight vector from the uniform
distribution on
\[
\left[
      - n^4,
      n^4
      \right]^{
(L-2) \cdot (k_0^2+k_0) + k_0 (d+2)+1
  }
,
\]
then condition (\ref{pth1eq1})
is satisfied for a weight vector $\bar{c}$ corresponding
to $X_1$ with probability at least
\[
\left(
\frac{1}{n^5}
\right)^{
(L-2) \cdot (k_0^2 + k_0) + k_0 \cdot (d+2) +1
}
     =
     \frac{1}{n^{
5 \cdot (L-2) \cdot (k_0^2 + k_0) + 5 \cdot k_0 \cdot (d+2) +5
     } }
     =: \eta_n
     .
\]
Hence after $
r_n = n \cdot \lceil \frac{1}{\eta_n} \rceil
$ of such independent choices (\ref{pth1eq1}) is never satisfied
with probability less than or equal to
\[
\left(
1-
\eta_n
\right)^{
r_n
}
\leq
\left(
1
-
\frac{n}{r_n}
\right)^{r_n}
\leq
\exp\left(
-
\frac{n}{r_n} \cdot r_n
\right)
=
e^{-n}.
\]
Now we consider $n$--times successively $r_n$ choices
of the weights, i.e.,
\[
k_n
=
n^2 \cdot
\lceil \frac{1}{\eta_n} \rceil
=
n^{5 \cdot (L-2) \cdot (k_0^2+k_0) + 5 \cdot k_0 \cdot (d+2)+7}
\]
such choices. Then the probability that in the first series
of weights there are no weights corresponding to $X_1$ chosen,
or in the second no weights corresponding to $X_2$, ...,
or in the $n$-th no weights corresponding to $X_n$ is bounded
from above
by
\[
\sum_{i=1}^n e^{-n} = n \cdot e^{-n}.
\]

Set
\[
L_n = n^{8 \cdot (L-2) \cdot (k_0^2+k_0) + 8 \cdot k_0 \cdot (d+2) + 16 \cdot L +15}.
\]
In the {\it third step of the proof} we show
that we have for $n$ sufficiently large
\begin{equation}
  \label{pth1eq5}
  \|\bc^{(t)}\|_\infty \leq 2 \cdot n^4
  \quad \mbox{for } t=0,1, \dots, t_n
\end{equation}
and
\begin{equation}
  \label{pth1eq6}
\| (\nabla_\bc F_n)(\bc) -  (\nabla_\bc F_n)(\bc^{(t)}) \|
\leq L_n \cdot \|\bc-\bc^{(t)}\|
\end{equation}
for all $\bc=\bc^{(t)}+s \cdot (\bc^{(t+1)}-\bc^{(t)})$
and all $s \in [0,1]$, for all $t=0,1, \dots, t_n-1$.

By Lemma \ref{le6} we know that
for $n$ sufficiently large
(\ref{le4eq1}) and (\ref{le4eq2})
hold for $c_3=1$ and $c_4=4$. The initial choice of our weights
implies furthermore (\ref{le4eq3}) and (\ref{le4eq4})
for $n$ sufficiently large. Application
of Lemma \ref{le4} yields (\ref{pth1eq5}). And (\ref{pth1eq5})
together with another application of Lemma \ref{le6} implies
(\ref{pth1eq6}).

In the {\it fourth step of the proof} we show
for $n$ sufficiently large
\begin{equation}
  \label{pth1eq7}
  F_n(\bc^{(t)}) \leq n^4
  \quad \mbox{for } t=0,1, \dots, t_n.
\end{equation}
Because of (\ref{pth1eq6}) we can conclude from Lemma \ref{le1}
that we have for $n$ sufficiently large
\[
 F_n(\bc^{(t+1)}) \leq  F_n(\bc^{(t)})  \quad \mbox{for } t=0,1, \dots, t_n-1.
\]
But the initial choice of the weights implies
\[
F_n(\bc^{(0)})
=
\frac{1}{n} \sum_{i=1}^n Y_i^2 \leq n^4.
\]

In the {\it fifth step of the proof} we show
that for $n$ sufficiently large and with probability
at least $1-n \cdot e^{-n}$ (\ref{le2eq1}) holds for all
$\bc=\bc^{(t)}$ $(t=0,1,\dots,t_n-1)$.
Because of the first and the second step of the proof it suffices
to show
\[
|\bar{c}_{j_i,k,l}^{(r)}-c_{j_i,k,l}^{(r)}| \leq
      \frac{1}{n}
\]
for all $i \in \{1, \dots, n\}$ and all $k,l,r$ with $r<L$,
where $\bar{c}_{j_i,k,l}^{(r)}$ and $c_{j_i,k,l}^{(r)}$ are the
corresponding components of $\bc^{(t)}$ and $\bc^{(0)}$.
Here $j_i$ is chosen such that
\begin{equation}
  \label{pth1eq8}
f_{j_i,j_i}^{(L)}(X_i) \geq 1 - e^{-n}
\quad \mbox{and} \quad
    \max_{t \in \{1, \dots,n\}, \, X_t \neq X_i} f_{j_i,j_i}^{(L)}(X_t)  \leq e^{-n}
.
\end{equation}
By Lemma \ref{le8} and the result of the fourth step of
the proof we can successively conclude
for $n$ sufficiently large that we have for
$t=0,1,\dots, t_n-1$
\[
\left|
\frac{\partial}{\partial_{j_i,k,l}^{(r)}} F_n(\bc^{(t)})
\right|
\leq
2 \cdot n^2 \cdot k_0^{L} \cdot (n^4)^{L+1} \cdot e^{-n}
\leq
\frac{1}{2 n^3}
=
\frac{1}{t_n \cdot \lambda_n} \cdot \frac{1}{n}
\]
and that consequently (\ref{pth1eq8}) holds for $\bc^{(t)}$.

In the {\it sixth step of the proof} we show
the assertion of Theorem \ref{th1}. By the results
of the third and the fifth step of the proof we know that
the assumptions of Lemma \ref{le3} are satisfied.
Application of Lemma \ref{le3} yields
\begin{eqnarray*}
  &&
F_n(\bc^{(t_n)})
-
\min_{g:\Rd \rightarrow \R}
\frac{1}{n} \sum_{i=1}^n |g(X_i)-Y_i|^2
\\
&&
\leq
  \left(
1-\frac{1}{2 \cdot n \cdot L_n}
\right)^{t_n}
\cdot
\left(
F_n(\bc^{(0)})
-
\min_{g:\Rd \rightarrow \R}
\frac{1}{n} \sum_{i=1}^n |g(X_i)-Y_i|^2
\right)
\\
&&
\leq
\exp \left(
- \frac{t_n}{2 \cdot n \cdot L_n}
\right)
\cdot
\left(
F_n(\bc^{(0)})
-
\min_{g:\Rd \rightarrow \R}
\frac{1}{n} \sum_{i=1}^n |g(X_i)-Y_i|^2
\right)
\\
&&
=
\exp(-n) \cdot
\left(
F_n(\bc^{(0)})
-
\min_{g:\Rd \rightarrow \R}
\frac{1}{n} \sum_{i=1}^n |g(X_i)-Y_i|^2
\right)
.
\end{eqnarray*}
With
\[
F_n(\bc^{(0)})
\leq
n^4
\]
we get the assertion. \hfill $\Box$

\subsection{Proof of Theorem \ref{th2}}
\label{se4sub2}

\begin{lemma}
\label{le9}
Let $n \in \N$, $(x_1,y_1)$, \dots, $(x_n,y_n) \in \Rd \times \R$,
$f:\Rd \rightarrow \R$, $\kappa_n>0$ and assume
\begin{equation}
\label{le9eq1}
\frac{1}{n} \sum_{i=1}^n | f(x_i)-y_i|^2
\leq
\min_{g:\Rd \rightarrow \R}
\frac{1}{n} \sum_{i=1}^n | g(x_i)-y_i|^2
+
\kappa_n.
\end{equation}
Set
\[
\bar{m}_n(x)
=
\frac{
\sum_{i=1}^n y_i \cdot I_{\{x_i=x\}}
}{
\sum_{i=1}^n I_{\{x_i=x\}}
}
\quad (x \in \Rd),
\]
where we use the convention $0/0=0$. Then we have for any $i \in \{1,
\dots,n\}$
\[
|f(x_i)-\bar{m}_n(x_i)|
\leq
\sqrt{ n \cdot \kappa_n}.
\]
\end{lemma}

\noindent
{\bf Proof.} We have
\begin{eqnarray*}
&&
\frac{1}{n} \sum_{i=1}^n | f(x_i)-y_i|^2
=
\frac{1}{n} \sum_{i=1}^n | f(x_i)-\bar{m}_n(x_i)|^2
+
\frac{1}{n} \sum_{i=1}^n | \bar{m}_n(x_i)-y_i|^2,
\end{eqnarray*}
since
\begin{eqnarray*}
&&
\frac{1}{n} \sum_{i=1}^n (f(x_i)-\bar{m}_n(x_i))\cdot
(\bar{m}_n(x_i)-y_i)
\\
&&
=
\frac{1}{n} \sum_{x \in \{x_1, \dots, x_n\}}
(f(x)-\bar{m}_n(x)) \cdot
\sum_{1 \leq i \leq n : x_i=x}
(\bar{m}_n(x_i)-y_i)
=
0.
\end{eqnarray*}
Application of (\ref{le9eq1}) yields
\[
\frac{1}{n} \sum_{i=1}^n | f(x_i)-\bar{m}_n(x_i)|^2
\leq \kappa_n,
\]
which implies the assertion.
\hfill $\Box$

\noindent
{\bf Proof of Theorem \ref{th2}.}
Set
\[
p_k=\frac{1}{n} \quad (k \in \{1, \dots, n\})
\]
and set $p_k=0$ for $k>n$. Set $x_k=(k/n,0,\dots,0)^T$ and define
the distribution of $(X,Y)$ by
\begin{enumerate}
\item
$\PROB[X=x_k]=p_k$ $(k \in \N)$,

\item
$Y=m(X)+\epsilon$ where $X$, $\epsilon$ are independent and
$m:\Rd \rightarrow \R$,

\item $\PROB\{\epsilon=-1\}=\frac{1}{2} = \PROB\{\epsilon=1\}$,

\item
$m(x)=0$ $(x \in \Rd)$.
\end{enumerate}
Then $m$ is the regression function of $(X,Y)$ and the distribution
of $(X,Y)$ satisfies the assumptions of Theorem \ref{th2}.

Set
\[
\bar{m}_n(x)
=
\frac{
\sum_{i=1}^n Y_i \cdot I_{\{X_i=x\}}
}{
\sum_{i=1}^n I_{\{X_i=x\}}
}
\quad (x \in \Rd).
\]
Using 
\[
| \bar{m}_n(x)|^2 \leq 2 \cdot |m_n(x_k)|^2 + 2 \cdot |m_n(x_k)-\bar{m}_n(x_k)|^2
\]
together with Lemma \ref{le9} we get
\begin{eqnarray*}
&&
\EXP \int | m_n(x)-m(x)|^2 \PROB_X (dx)
\\
&&
\geq
\EXP \left\{
\sum_{k=1}^{n}
|m_n(x_k)|^2 \cdot p_k \cdot I_{\{\sum_{i=1}^n I_{\{X_i=x_k\}}>0\}} \cdot I_{\{ U \in C_n\}}
\right\}
\\
&&
\geq
\EXP \left\{
\sum_{k=1}^{n}
\left(
\frac{1}{2} |\bar{m}_n(x_k)|^2
-|m_n(x_k) - \bar{m}_n(x_k)|^2
\right)
\cdot p_k \cdot  I_{\{\sum_{i=1}^n I_{\{X_i=x_k\}}>0\}}  \cdot I_{\{ U \in C_n\}}
\right\}
\\
&&
=
\EXP \left\{
\sum_{k=1}^{n}
\frac{1}{2} |\bar{m}_n(x_k)|^2
\cdot p_k \cdot I_{\{\sum_{i=1}^n I_{\{X_i=x_k\}}>0\}} 
\right\}
\\
&&
\hspace*{1cm}
-
\EXP \left\{
\sum_{k=1}^{n}
\frac{1}{2} |\bar{m}_n(x_k)|^2
\cdot p_k \cdot I_{\{\sum_{i=1}^n I_{\{X_i=x_k\}}>0\}}   \cdot I_{\{ U \in C_n^c\}}
\right\}
\\
&&
\hspace*{1cm}
-
\EXP \left\{
\sum_{k=1}^{n}
|m_n(x_k) - \bar{m}_n(x_k)|^2
\cdot p_k \cdot I_{\{\sum_{i=1}^n I_{\{X_i=x_k\}}>0\}}  \cdot I_{\{ U \in C_n\}}
\right\}
\\
&&
\geq
\frac{1}{2}
\cdot
\sum_{k=1}^{n}
\EXP \left\{
|\bar{m}_n(x_k)|^2
\cdot I_{\{\sum_{i=1}^n I_{\{X_i=x_k\}}>0\}} 
\right\}
\cdot p_k
- \frac{1}{2} \cdot \PROB_U(C_n^c) - n \cdot \kappa_n.
\end{eqnarray*}
The definition of $\bar{m}_n$ implies
\begin{eqnarray*}
&&
\sum_{k=1}^{n}
\EXP \left\{
|\bar{m}_n(x_k)|^2
\cdot I_{\{\sum_{i=1}^n I_{\{X_i=x_k\}}>0\}} 
\right\}
\cdot p_k
\\
&&
\geq
\sum_{k=1}^{n}
\EXP \left\{\EXP \left\{
|\bar{m}_n(x_k)|^2 \big| X_1, \dots, X_n \right\}
\cdot I_{ \{ \sum_{i=1}^n I_{\{X_i=x_k\}}>0 \}}
\right\}
\cdot p_k
\\
&&
=
\sum_{k=1}^{n}
\EXP \left\{
\frac{1}{
\sum_{i=1}^n I_{\{X_i=x_k\}}
}
\cdot I_{ \{ \sum_{i=1}^n I_{\{X_i=x_k\}}>0 \}}
\right\}
\cdot p_k
.
\end{eqnarray*}
Using the fact that $\sum_{i=1}^n I_{\{X_i=x_k\}}$ is binomially
distributed
with $n$ degrees of freedom and probability of success $p_k$ we get
\begin{eqnarray*}
&&
\sum_{k=1}^{n}
\EXP \left\{
\frac{1}{
\sum_{i=1}^n I_{\{X_i=x_k\}}
}
\cdot I_{ \{ \sum_{i=1}^n I_{\{X_i=x_k\}}>0 \}}
\right\}
\cdot p_k
\\
&&
=
\sum_{k=1}^{n}
\sum_{i=1}^n
\frac{1}{i} \cdot \left( {n \atop i} \right) p_k^i \cdot (1-p_k)^{n-i}
\cdot p_k
\\
&&
\geq
\sum_{k=1}^{n}
\sum_{i=1}^n
\frac{1}{i+1} \cdot \left( {n \atop i} \right) p_k^i \cdot (1-p_k)^{n-i}
\cdot p_k
\\
&&
=
\frac{1}{n+1} \cdot
\sum_{k=1}^{n}
\sum_{i=1}^n
\left( {n+1 \atop i+1} \right) p_k^{i+1} \cdot (1-p_k)^{n+1-(i+1)}
\\
&&
=
\frac{n}{n+1} 
\cdot
\left(1-\left(1-\frac{1}{n} \right)^{n+1} - (n+1) \cdot \frac{1}{n} \cdot
\left(1-\frac{1}{n} \right)^{n}\right)
\\
&&
\geq
\frac{n}{n+1}
\cdot
\left(1-\frac{2n+1}{n} \cdot
\left(1-\frac{1}{n} \right)^{n}\right)
\\
&&
\geq
\frac{10}{11}
\cdot
\left(1-\frac{21}{10} \cdot
\frac{1}{e} \right),
\end{eqnarray*}
where the last inequality holds for $n \geq 10$.

Putting together the above results implies the assertion.
\hfill $\Box$

\end{document}